\documentclass{article}

\usepackage{arxiv}

\usepackage[utf8]{inputenc} % allow utf-8 input
\usepackage[T1]{fontenc}    % use 8-bit T1 fonts
\usepackage{hyperref}       % hyperlinks
\usepackage{url}            % simple URL typesetting
\usepackage{booktabs}       % professional-quality tables
\usepackage{amsfonts}       % blackboard math symbols
\usepackage{nicefrac}       % compact symbols for 1/2, etc.
\usepackage{microtype}      % microtypography
\usepackage{lipsum}
\usepackage{graphicx}
\graphicspath{ {./images/} }

\title{On the {G}aussian Approximation to {B}ayesian Posterior Distributions}

\author{
 Christoph Fuhrmann \\
  School of Education\\
  Institute for educational research\\
  University of Wuppertal\\
  D-42097 Wuppertal \\
  \texttt{christoph.fuhrmann@uni-wuppertal.de} \\
   \And
  Hanns Ludwig Harney \\Max-Planck-Institute for Nuclear Physics\\ Postfach 103980\\
  D-69029 Heidelberg\\
  \texttt{hanns-ludwig.harney@mpi-hd.mpg.de} \\ 
   \And
  Klaus Harney \\
  Faculty of Philosophy and Educational Science\\
  Institute for educational research \\   Ruhr-Universit\"at \\ D-44801 Bochum \\
   \texttt{klaus.harney@ruhr-uni-bochum.de} \\  
  \And
   Andreas M{\"u}ller \\
  Faculty of Sciences/Physics Department and\\ 
  Institut Universitaire de Formation des Enseignants\\ University of Geneva\\CH - 1211 Genève 4 \\
  \texttt{andreas.mueller@unige.ch} \\
}

\renewcommand{\undertitle}{}
\renewcommand{\headeright}{}

\begin{document}
\maketitle
\begin{abstract}
  The present article derives the minimal number $N$ of
  observations needed to consider a Bayesian posterior distribution as
  Gaussian. Two examples are presented. Within one of them,
  a chi-squared distribution, the observable $x$ as well as the parameter
  $\xi$ are defined all over the real axis, in the other one, the binomial distribution, the
  observable $x$ is an entire number while the parameter $\xi$ is defined
  on a finite interval of the real axis.
  The required minimal $N$ is high in the first case and
  low for the binomial model. In both cases the precise definition of the
  measure $\mu$ on the scale of $\xi$ is crucial.
\end{abstract}

% keywords can be removed
\keywords{{B}ayesian posterior \and {G}aussian approximation \and chi-squared and binomial
  distributions}

\section{General Notions and Definitions}
\label{sec:1}

Bayesian statistics distinguishes the observations $x_1\, ,\dots\, ,x_N\, $
from the parameter $\xi$ that ``conditions'' them. 
{B}ayes' theorem --- given below --- expresses the uncertainty about $\xi$
via a probability distribution of $\xi$ conditioned by the observations
$x_1\, ,\dots\, ,x_N\, .$ This so-called posterior distribution becomes always narrower with increasing $N\, .$ 
By consequence the ``true value'' of $\xi$ is approached more and more closely. Simultaneously the
posterior approaches a Gaussian distribution. This is a consequence of the fact that
the posterior after $N$ observations becomes the $N$-th power of the posterior from one observation.

The present article, based on Bayesian statistics, derives the minimal $N$ needed for the
Gaussian approximation.
In the first part, comprising Sects. \ref{sec:2} and \ref{sec:3}, the variables $x$
and $\xi$ are defined along the real axis as a whole. A chi-squared distribution provides an
example. In the second part, i.e. in Sects. \ref{sec:4} and \ref{sec:5}, the
parameter $\xi$ is defined on a finite interval of the real axis. The
so-called trigonometric distribution serves as example.

We explain some notions used throughout the following text.

A {\it statistical model} $p(x|\xi )$ formulates the relation between the
observations $x$ and the parameter $\xi\, .$ 
A statistical model is normalised according to
\begin{equation}
  \int_{domain\, x}{\rm d}x\, p(x|\xi )\, =\, 1
  \label{1.0}
\end{equation}
for every $\xi\, .$ Here, ``$domain\, x$'' means that the integral extends over the domain of definition of the observed $x\, .$

{\it {B}ayes' theorem} \cite{Bayes:1763} states the posterior distribution
\begin{equation}
  P(\xi |x)\, =\, \frac{p(x|\xi )\mu (\xi )}{m(x)}\, ,
  \label{1.1}
\end{equation}
of $\xi$ conditioned by $x\, .$ The Bayesian prior distribution $\mu (\xi )$ serves also as the measure of integration over $\xi\, ,$
see Eq. (\ref{4.4a}) and the explanation there. The quantity
\begin{equation}
  m(x)\, =\, \int_{domain\, \xi}{\rm d}\xi\, p(x|\xi )\mu (\xi )
  \label{1.2}
\end{equation}
normalises $P(x|\xi )$ to unity.
Equation (\ref{1.1}) gives the posterior in the case of a single
observation $x\, .$ There is usually more than one
observation, but in the present text there shall be only one parameter $\xi\, .$
For the case of $N$ observations $\bold{x}=(x_1,\dots ,x_N)$ the posterior
is formulated in Sects. \ref{sec:2.4} and \ref{sec:4.4}.

The {\it Fisher information} \cite{Fisher:22,Fisher:25,Fisher:59,Fisher:67}
is the expectation value
\begin{equation}
  F(\xi )\, =\, \int_{domain\, x}{\rm d}x\, p(x|\xi )
  \left[\frac{\partial}{\partial\xi}\ln p(x|\xi )\right]^2
  \label{1.3}
\end{equation}
of the quantity $[\partial/\partial\xi\, \ln p(x|\xi )]^2\, .$ The Fisher information is always positive. 
When $p$ is a Gaussian distribution then $F$ is the inverse of the mean
square value of the Gaussian.

{\it Form invariance} is a symmetry relation between the observation $x$ and
the condition $\xi\, .$ It is defined by a group --- in the sense of the
theory of Lie groups \cite{Haar:33,Jeffreys:46} --- of transformations that
leave $p(x|\xi )$ unchanged when applied to both, $x$ and $\xi\, .$ See Chap. 6 of \cite{Harney:16}
and the textbooks \cite{Adams:69,Broecker:85}.

In the examples of the present text the
symmetry shows up by the fact that the model $p(x|\xi )$ depends on the difference $x-\xi\, .$
The symmetry group then consists of all possible
simultaneous translations of $x$ and $\xi$ by the same shift.
Many statistical models can be reformulated such as to display translational form invariance.
In section \ref{sec:2} the details of translational form invariance are discussed.
It is taken as the starting point for a Gaussian approximation to the posterior distribution
since the Gaussian itself is form invariant under translations. The invariant measure of the group
of translations is identified with the prior distribution. Here, this implies
\begin{equation}
  \mu (\xi )\, \equiv\, {\rm const}\, .
  \label{1.5}
\end{equation}
In Sects. \ref{sec:2} and \ref{sec:4} two different cases are considered: Continuous $x$ and $\xi$ within
a chi-squared model and dichotomic $x$ within the binomial model. Both cases are widespread and of practical
interest for many applications.

A {\it likelihood function} ${\cal L}_N(\xi )$ is proportional to the probability density $p_N(\bold{x}|\xi )$ of $N$ observations $\bold{x}$
considered as a function of $\xi$ while $\bold{x}$ is given. 
When --- in the present context of Eq. (\ref{1.5}) --- the domain of definition of
$\xi$ extends over the whole real axis, the likelihood function possesses a maximum: Since the posterior
is normalised, it must tend to zero when $\xi$ goes to infinity. 

The value $\xi^{\rm ML}$ where the maximum occurs, is called the
{\it maximum likelihood} or ML {\it estimator} of the ``true value'' of $\xi\, .$
For every series of observations $\bold{x}$ there is a ML estimator
$\xi^{\rm ML}=\xi^{\rm ML}(\bold{x})\, .$ In sections \ref{sec:2.4} and \ref{sec:4.4} one shall see that $\xi^{\rm ML}$ is the
{\it sufficient statistic} \cite{Kendall:76} of the model --- a notion widely
discussed in the development of the Rasch model
\cite{Rasch:60,Rasch:66a,Rost:04,Fisch:81,Fisch:95,Fuhrm:16,Harney:16}.

It might appear that the constancy of the prior is nothing but the
"principle of indifference", well-known in the history of statistical and
Bayesian reasoning, well-known also for its difficulties when applied to
transform probability densities \cite{Edwards:92,Maxent:89}. Note that in the present paper the constancy
of the prior occurs as a consequence of the well-defined group theoretical property of ``form invariance'' \cite{Hartigan:64}.

In Sects. \ref{sec:3.1} and \ref{sec:5} the likelihood function of a model
$p(x|\xi )$ is compared to the Gaussian model
\begin{equation}
  G(x|\xi )\, =\, (2\pi\sigma^2)^{-1/2}\exp\left(-\frac{(x-\xi )^2}{2\sigma^2}\right)\, ,\quad\quad -\infty <x,\xi <\infty\, .
  \label{1.6a}
\end{equation}
The $N$-fold Gaussian model, i.e. the distribution of $N$ observations $\bold{x}=(x_1,\dots ,x_N)\, ,$ is
\begin{equation}
  G_N(\bold{x}|\xi )\,=\, \left(\frac{1}{2\pi\sigma^2}\right)^{N/2}\prod_{k=1}^N\, \exp\left(-\frac{(x_k-\xi )^2}{2\sigma^2}\right)\, .
\label{1.6b}
\end{equation}
This yields the posterior distribution
\begin{equation}
  G_N(\xi |\bold{x})\, =\, \left(\frac{N}{2\pi\sigma^2}\right)^{1/2}\exp\left(-\frac{N}{2\sigma^2}(<x>-\xi )^2\right)\, ,
  \label{1.6}
\end{equation}
where $<x>$ is the average
\begin{equation}
  <x>\, =\, \frac{1}{N}\sum_{k=1}^N\, x_k\, ,
  \label{1.6c}
\end{equation}
see appendix \ref{I}. The Fisher information of the Gaussian (\ref{1.6a}) is
\begin{equation}
  F^{\rm Gauss}\, \equiv\, \sigma^{-2}\, .
  \label{1.6d}
\end{equation}

With increasing number $N$ of observations 
the posterior of any model $p$ with the prior (\ref{1.5}) assumes the
Gaussian form and becomes always narrower tending towards Dirac\rq s delta distribution.
We shall determine the minimal $N$ which allows to approximate a given posterior by the Gaussian (\ref{1.6})
within the interval $|\xi^{\rm ML}-\xi|<3\sigma /\sqrt{N}\, .$ This interval
contains $99.73$ percent of the area under the Gaussian function.
Two examples are studied, the chi-squared model in Sect. \ref{sec:3}
and the binomial model in Sect. \ref{sec:5}. In the two cases the minimal $N$ strongly differ.

\section{Form Invariance Along the Real Axis}
\label{sec:2}
The present section considers a statistical model $p(x|\xi )$ where
both, the observable $x$ and the parameter $\xi\, ,$ are defined all along the
real axis. Furthermore, $x$ and $\xi$ shall be related via form
invariance. These two properties allow a general approximation to the logarithmic likelihood function
which contains the sum over $\ln p(x_k-\xi )\, .$ These logarithms are $N$ random numbers since the $x_k$
are random, whence the sum over the $\ln p(x_k-\xi )$ will have a Gaussian distribution for sufficiently large $N\, .$
Yet the central limit theorem does not allow to answer the question, how large $N$ must be for the Gaussian
approximation to apply. However, the sum over the $\ln p(x_k-\xi )$ becomes equal to the $N$-fold (negative) Kullback-Leibler divergence
$H(\xi^{\rm ML}|\xi )\, ,$ see Ref. \cite{Kullback:51} and Eq. (\ref{3.4}). The Kullback-Leibler divergence is a quantity that measures the distance
between the distributions $p(x|\xi^{\rm ML})$ and $p(x|\xi )\, .$
It will be represented by a Taylor expansion with respect to $\xi\, .$
With increasing $N$ its terms of higher order become negligible as compared to the term of second order.
This leads to the criterion for the Gaussian approximation.

\subsection{The prior distribution}
\label{sec:2.1}
The above-mentioned two properties of the model $p(x|\xi )$ entail that
$p$ depends on the difference $x-\xi$ and only on this difference.
Thus the model reads
\begin{equation}
  p(x|\xi )\, =\, p(x-\xi )\, ,\quad\quad -\infty <x,\xi <\infty\, .
  \label{2.1}
\end{equation}

Neither Bayes \cite{Bayes:1763} nor Laplace \cite{Laplace:1774} who
independently established Bayes' theorem, gave a prescription to obtain the
prior distribution. Form invariance gives us the prescription: The prior shall
be invariant under the symmetry group of translations
\cite{Jeffreys:46,Kass:90,Kass:96,Harney:16}.

Under the translation
\begin{equation}
  \xi^{\prime}\,  =\, \xi +a
  \label{2.2}
\end{equation}
the prior transforms as a density, i.e. $\mu (\xi )$ goes over into
\begin{eqnarray}
  \mu_{T}\, &=&\, \mu(\xi^{\prime}-a)\left|\frac{{\rm d}\xi}{{\rm d}\xi^{\prime}}
  \right|\nonumber\\
  &=&\, \mu (\xi^{\prime}-a)\, .
  \label{2.3}
\end{eqnarray}
This shall be independent of $a\, ;$ thus $\mu$ is constant as foreseen in
Eq. (\ref{1.5}). See also Eq. (\ref{4.4a}) in Sect. \ref{sec:2.3}.

For $N$ events $\bold{x}\, ,$ conditioned by one and
the same value of $\xi\, ,$ the model is
\begin{equation}
  p_N(\bold{x}|\xi )\, =\, \prod_{k=1}^N\, p(x_k-\xi )\, ,\quad\quad
  -\infty <x_k,\xi <\infty\, .
  \label{2.4}
\end{equation}
The posterior is
\begin{equation}
  P_N(\xi |\bold{x})\, =\, \mu (\xi )\frac{\prod_{k=1}^N\, p(x_k-\xi )}
  {m(\bold{x})}\, ,
  \label{2.5}
\end{equation}
where
\begin{equation}
  m(\bold{x})\, =\, \int_{-\infty}^{\infty}{\rm d}\xi\, \mu (\xi )
  \prod_{k=1}^N\, p(x_k-\xi )\, .
  \label{2.6}
\end{equation}
Since $\mu$ is constant its value drops out of the posterior and $P$ becomes
\begin{equation}
  P_N(\xi |\bold{x})\, =\, \frac{\prod_{k=1}^N\, p(x_k-\xi )}
  {\int_{-\infty}^{\infty}{\rm d}\xi^{\prime}\, \prod_{k=1}^N\, p(x_k-\xi^{\prime})}\, .
  \label{2.7}
\end{equation}

One can numerically calculate this expression, determine the shortest interval
in which $\xi$ is found with the probability of $99.73$ percent, and thus obtain an error
interval for $\xi\, ,$ see Chap. 3 of \cite{Harney:16}.
In Sect. \ref{sec:2.4} this procedure is replaced by considering the logarithm of the likelihood function.
This will show that
the posterior $P_N$ of Eq. (\ref{2.7}) tends towards a Gaussian with increasing $N\, .$

\subsection{The Maximum-Likelihood Estimator}
\label{sec:2.2}

The expression $p_N(\bold{x}|\xi )$ of Eq. (\ref{2.4}) is called a
likelihood function ${\cal L}_N(\xi )$ when it is considered as a function
of $\xi\, ,$ while $\bold{x}$ is given.

A likelihood function in the context of translational invariance possesses a
maximum since it is a normalisable function defined all along the real axis.
When there are several maxima, one must look for the absolute maximum or even
redefine the model such that there is an absolute maximum. The place
$\xi^{\rm ML}\, ,$ where the maximum occurs, depends on the observed
$\bold{x}\, .$ Hence, $\xi^{\rm ML}=\xi^{\rm ML}(\bold{x})$ is a function of
$\bold{x}\, .$ The value $\xi^{\rm ML}$ is called the maximum likelihood (ML)
estimator of the parameter $\xi\, .$ For the example of the chi-squared model
in Sect. \ref{sec:3.2} the ML estimator is given in appendix
\ref{B}. 

The ML estimator has been introduced by R.A. Fisher
\cite{Fisher:12,Fisher:22,Aldrich:97}. It estimates the ``true value'' of
$\xi$ which conditions the observations $\bold{x}\, .$ For every finite $N\, ,$
however, the true value remains hidden. With $N\to\infty$ the ML estimator
converges to it. An example is given by the Gaussian model (\ref{1.6b}); the
average (\ref{1.6c}) is its ML estimator.

Section \ref{sec:2.4} shows that $\xi^{\rm ML}(\bold{x})$ is the sufficient
statistic of the model $p\, ,$ i.e. the observations $\bold{x}$ enter into
the posterior distribution only via $\xi^{\rm ML}\, .$ Compare page 22 of
\cite{Kendall:51} and Sect. 3.1.3 of \cite{Fuhrm:16} and chapters 2 and 3 of
Ref. \cite{Harney:16}.

Neyman and Scott \cite{NeymanScott:48} have argued against ML estimation. Their argument says that a bias may remain between the expectation value of $\xi$ and the ML estimator. This has caused a considerable debate \cite{Stigler:07,Spanos:13}. The argument of Neyman and Scott is bound to  the distinction between between a ``structural'' and an ``incidental'' parameter. The model they studied, is form invariant which means that there is a Lie group of transformations that leaves it invariant. Each of the two parameters describes a subgroup. The elements of the different subgroups do not commute with each other. One of the subgroups describes translations, the other one describes dilations, see, e.g., Sect. 7.3 of \cite{Harney:16}. The present article describes a more basic situation. We also have two parameters; however, the present subgroups are identical. Both are translational.

At least in the present context of translational form invariance such a bias
becomes arbitrarily small with increasing $N$ --- as is shown for the chi-squared model of Eq. (\ref{6.11})
in Sect. \ref{sec:3.2}. It possesses such a bias; but the fact that
the posterior tends to a Gaussian with increasing $N$ implies that the bias goes to zero. Future research sould show
whether this holds also for the model discussed by Neyman and Scott.

We study the transition  of the posterior to a Gaussian distribution by help of the logarithm of $P_N$ 
given in Eq. (\ref{2.5}). Since both, $\mu (\xi )$ and $m(\bold{x})\, ,$ are independent of $\xi\, ,$
this likelihood function is proportional to the product of the $p(x_k-\xi )\, ,$ i.e.
\begin{equation}
  {\cal L}_N(\xi )\, \propto\, \prod_{k=1}^N\, p(x_k-\xi )\, .
  \label{3.0}
\end{equation}
Therefore the logarithm of ${\cal L}_N$ is --- up to an additive constant ---
the sum over the logarithms $\ln p(x_k-\xi )\, ,$
\begin{equation}
  {\rm const}+\ln{\cal L}_N(\xi )\, =\, \sum_{k=1}^N\, \ln p(x_k-\xi )\, .
  \label{3.0a}
\end{equation}
For sufficiently large $N$ the sum over the logarithms is expressed by the expectation value
of $\ln p(x-\xi )$ which in principle requires
\begin{equation}
  {\rm const}+\ln{\cal L}_N(\xi )\, =\, N\int_{\infty}^{\infty}{\rm d}x\, p(x-\xi^{\rm true})\ln p(x-\xi)\, .
  \label{3.0bb}
\end{equation}
This integral is the expectation value of $\ln p(x-\xi )$ taken with the distribution conditioned by the ``true value'' of $\xi\, .$
Jaynes and Bretthorst \cite{Jaynes:12} have called it the ``asymptotic'' likelihood function,
asymptotic in the sense of $N\to\infty\, .$
The true value of $\xi$ remains, however, hidden for every finite $N\, .$ We replace it by the ML estimator $\xi^{\rm ML}$
obtained from the observations $\bold{x}=(x_1,\dots ,x_N)\, .$ Then the ``asymptotic'' likelihood function becomes
\begin{equation}
  {\rm const}+\ln{\cal L}_N(\xi )\, =\, N\int_{\infty}^{\infty}{\rm d}x\, p(x-\xi^{\rm ML})\ln p(x-\xi)\, .
  \label{3.1aa}
\end{equation}
We want to define $\rm const$ such that $\ln{\cal L}_N$ becomes a Kullback-Leibler distance. This is reached
when $\rm const$ is set to
\begin{equation}
  {\rm const}\, =\, -N\int_{\infty}^{\infty}{\rm d}x\, p(x-\xi^{\rm ML})\ln p(x-\xi^{\rm ML})
  \label{3.0b}
\end{equation}
Then the ``asymptotic'' likelihood function becomes
\begin{equation}
  \ln {\cal L}_N(\xi )\, =\, N\int_{{\rm domain}\, x}{\rm d}x\, p(x-\xi^{\rm ML})\ln\frac{p(x-\xi )}{p(x-\xi^{\rm ML})}\, .
  \label{3.1a}
\end{equation}
We call this integral the functional
\begin{equation}
  H(\xi^{\rm ML}|\xi )\, =\, \int_{{\rm domain}\, x}{\rm d}x\, p(x-\xi^{\rm ML})\ln\frac{p(x-\xi )}{p(x-\xi^{\rm ML})} \, .
  \label{3.4}
\end{equation}
It is the negative of the Kullback-Leibler divergence \cite{Kullback:51,Hobson:73,Karmeshu:03} between the distributions conditioned by $\xi^{\rm ML}$ and by $\xi\, .$
Thus the logarithmic likelihood function becomes 
\begin{equation}
  \ln{\cal L}_N(\xi )\, =\, N\, H(\xi^{\rm ML}|\xi )\, .
  \label{3.2}
\end{equation}
The Taylor expansion of $H$ with respect to $\xi$ will yield our criterion by requiring that terms of higher than the second order
be negligible.

How large must $N$ be for $\ln{\cal L}_N$ to become equal to the expectation value (\ref{3.1a})?
We follow the assumption that there is such an $N\, .$ For this and larger $N$ equation (\ref{3.1a}) holds: The
logarithmic likelihood function of the Gaussian (\ref{1.6}) is
\begin{equation}
  \ln{\cal L}_N^{\rm Gauss}(\xi )\, =\, N\left[-\frac{(<x>-\xi )^2}{2\sigma^2}\right]\, ,
  \label{1.6e}
\end{equation}
where the quantity in rectangular brackets is the functional (\ref{3.4}) for the Gaussian (\ref{1.6a}),
\begin{equation}
  H^{\rm Gauss}(\xi^{\rm ML}|\xi )\, =\, -\frac{(<x>-\xi )^2}{2\sigma^2}\, .
  \label{1.6f}
\end{equation}
The quantity $<x>\, ,$ given in Eq. (\ref{1.6c}), is the ML estimator of the Gaussian model.
One sees that Eqs.\,  (\ref{3.1aa}) and (\ref{3.2}) are fulfilled by the Gaussian model.

Equation (\ref{A.7}) in appendix \ref{A} shows that the maximum value of $H(\xi^{\rm ML}|\xi )$ --- and thus of the logarithmic likelihood ---
occurs at $\xi =\xi^{\rm ML}$  for every distribution of the form $p(x-\xi )\, .$

According to Eq. (\ref{3.2}) the likelihood ${\cal L}_N(\xi )$ has the form of the $N$-th power of ${\cal L}_{N=1}(\xi )$
since $\ln{\cal L}_N(\xi )$ is proprotional to $N\, .$ The translational invariance of $p$ given by Eq. (\ref{2.1}) leads to a translational
invariance of $H$ since
\begin{eqnarray}
  H(\xi^{\rm ML}|\xi )\, &=&\, \int_{-\infty}^{\infty}{\rm d}x\, p(x-\xi^{\rm ML})\ln\frac{p(x-\xi )}{p(x-\xi^{\rm ML})}\nonumber\\
  &=&\, \int_{-\infty}^{\infty}{\rm d}x^{\prime}\, p(x^{\prime}-\xi^{\rm ML}+\xi )\ln\frac{p(x^{\prime})}{p(x^{\prime}-\xi^{\rm ML}+\xi )}\nonumber\\
  &=&\, H(\xi^{\rm ML}-\xi )|0)\, .
  \label{3.4b}
\end{eqnarray}
This result has been obtained by substituting
\begin{equation}
  x^{\prime}\, =\, x-\xi\, .
  \label{3.4a}
\end{equation}
Thus $H(\xi^{\rm ML}|\xi )$ depends on the difference $\xi^{\rm ML}-\xi$ and only on this difference.

When there are isolated points where $p(x^{\prime})$ or $p(x^{\prime}-\xi^{\rm ML}+\xi )$ vanish, then appendix \ref{H} shows that the integral
exists despite the divergence of the integrand in (\ref{3.4b}).

\subsection{The Fisher Information}
\label{sec:2.3}

The information that carries the name of R.A. Fisher \cite{Fisher:25} is a central concept of statistics and estimation theory.
It provides the prior, if form invariance holds. The Fisher information $F$ is defined as the expectation value
of the squared derivative of the logarithmic likelihood function, see Refs. \cite{Fisher:22,Fisher:25,Fisher:59,Fisher:67} and Sect. 4.2 of
\cite{Fuhrm:16}, i.e.
\begin{eqnarray}
  F(\xi )\, &=&\, \int_{domain\, x}{\rm d}x\, p(x|\xi )
  \left[\frac{\partial}{\partial\xi}\ln p(x|\xi )\right]^2\nonumber\\
            &=&\, -\int_{domain\, x}{\rm d}x\,
  p(x|\xi )\frac{\partial^2}{\partial\xi^2}\ln p(x|\xi )\, .
  \label{4.1}
\end{eqnarray}
The two expressions equal each other because $p$ is normalised to unity for
every $\xi\, ,$ see appendix \ref{C}. Note that the definition (\ref{4.1}) refers to a single observation, see \cite{LC:98}.
Other authors, however, define the Fisher information for multiple observations \cite{PP:02}.
The integration in (\ref{4.1}) extends over the
full domain of definition of $x\, .$ Given the domain of Eq. (\ref{2.1}) one has by use of the substitution (\ref{3.4a})
\begin{eqnarray}
  F\, &=&\, -\int_{-\infty}^{\infty}{\rm d}x\,
  p(x-\xi )\frac{\partial^2}{\partial\xi^2}\ln p(x-\xi )\nonumber\\
  &=&\, -\int_{-\infty}^{\infty}{\rm d}x^{\prime}\, p(x^{\prime})
  \frac{{\rm d}^2}{{\rm d}{x^{\prime}}^2}\ln p(x^{\prime})\, .
  \label{4.2}
\end{eqnarray}
Thus the Fisher information of the model (\ref{2.1}) is independent of
$\xi\, .$ This is a consequence of form invariance under translations.

By the definitions of $H$ in (\ref{3.4}) and $F$ in (\ref{4.1}) one
recognises that
\begin{equation}
  F\, =\, -\left.\frac{\partial^2}{\partial\xi^2}H(\xi^{\rm ML}|\xi )
    \right|_{\xi =\xi^{\rm ML}}\, .
    \label{4.3}
\end{equation}
Since $F$ is positive, the second derivative of $H$ is seen to be
negative at $\xi =\xi^{\rm ML}\, .$ This means that the
second derivative of ${\cal L}(\xi )$ is negative, as it should be at a
maximum.

The prior distribution $\mu$ is proportional to the square root of $F\, ,$
\begin{equation}
  \mu\, \propto\, \sqrt{F}\, ,
  \label{4.4a}
\end{equation}
according to Hartigan \cite{Hartigan:64} and Jaynes \cite{Jaynes:68}.
A measure is the inverse unit of length on the scale of $\xi\, .$

By identifying the integration measure on the scale of $\xi$ with the Bayesian
prior distribution, we follow Jeffreys \cite{Jeffreys:46} as was done by
other authors \cite{Kass:80,Amari:85,Kass:96,Rao:85}.
The prior is independent of $N\, .$ Since it is
independent of $\xi$ and  $\xi^{\rm ML}\, ,$ a given
difference $\xi^{\rm ML}-\xi$ describes the same length of way everywhere on the scale
of $\xi\, .$ See also \cite{Aykac:77,BayStat3,Zellner:77,Zellner:96,Harney:16} and the footnote
\footnote{The prior distribution by Zellner \cite{Zellner:77} is constructed
  by help of Shannon's information. For the model (\ref{2.1}) Zellner agrees
  with the measures by Jeffreys and Kass in that Zellner's prior is constant.
  However, his prior does not behave as a density under transformations. So it
  is not generally proportional to the square root of (\ref{4.3}).}.

\subsection{The Posterior of a Form Invariant Model}
\label{sec:2.4}

>From Eq. (\ref{3.2}) and the fact that the prior is constant follows that the
posterior distribution is
\begin{equation}
  P_N(\xi |\bold{x})\, =\, {\cal N}_N\exp\Big(NH(\xi^{\rm ML}|\xi )\Big)\, ,
  \label{5.4}
\end{equation}
where ${\cal N}_N$ normalises $P_N$ to unity. 
The posterior from $N$ observations has the form of the
$N$-th power of the posterior for $N=1\, .$
Therefore with increasing $N$ the posterior
is more and more restricted to the immediate neighborhood of the maximum at $\xi =\xi^{\rm ML}\, .$
At this maximum $H(\xi^{\rm ML}|\xi^{\rm ML})$ vanishes, since $H$ equals the Kullback-Leibler distance between
the distributions $p(x-\xi^{\rm ML})$ and $p(x-\xi )\, .$ Thus for $\xi =\xi^{\rm ML}$ the exponential in
(\ref{5.4}) has the value of unity,
\begin{equation}
  \exp \left(NH(\xi^{\rm ML}|\xi^{\rm ML} )\right)\, =\, 1\, .
  \label{5.5}
\end{equation}
With increasing $N$ the curvature of the likelihood function $\exp (NH(\xi^{\rm ML}|\xi ))$  at its
maximum increases according to $NF\, ,$ were $F$ is the Fisher information at $N=1\, .$ The
Fisher information is independent of $\xi$ in the present context. Thus the likelihood function
behaves as the Gaussian (\ref{1.6}), see Eqs. (\ref{1.6e}) and (\ref{1.6f}). It does so in a
suitable interval around the maximum value (\ref{5.5}) of the likelihood. Thus within that interval
the likelihood becomes a Gaussian function.

One also sees that the posterior
depends on the observations $\bold{x}$ via and only via the estimator
$\xi^{\rm ML}=\xi^{\rm ML}(\bold{x})\, ;$ this means that $\xi^{\rm ML}(\bold{x})$
is the sufficient statistic \cite{Kendall:51}. 
The notion of sufficient statistic has been widely discussed in the development
of the Rasch model of item response theory
\cite{Rasch:60,Rasch:66a,Fisch:81,Fisch:95,Rost:04,Fuhrm:16,Harney:16}.
We shall come to it in Sect. \ref{sec:4.5} where the binomial model is treated.

\section{The Criterion for the {G}aussian Approximation}
\label{sec:3}

A criterion will now be formulated that allows to find the minimal $N$ so that
the Gaussian distribution (\ref{1.6})
can be accepted as approximating the posterior $P_N$ in Eq. (\ref{5.4}).

The Gaussian approximation is justified
whenever the Taylor expansion of $\ln{\cal L}_N(\xi )$ up to the second order in
$\xi$ is ``sufficiently precise''. This is the expansion 
of $H(\xi^{\rm ML}|\xi )$ with respect to $\xi$ at the point
$\xi =\xi^{\rm ML}\, .$ In Sect. \ref{sec:3.1} the Taylor
expansion is written down and our understanding of ``sufficiently precise''
is defined. This yields the desired criterion. In Sect. \ref{sec:3.2}, the
criterion is applied to a chi-squared model.

\subsection{Formulation of the Criterion}
\label{sec:3.1}

The Taylor expansion of $H$ up to an arbitrary order $n$ is
\begin{equation}
  H(\xi^{\rm ML}|\xi )\, =\, \sum_{\nu =1}^n\,
  \frac{(\xi -\xi^{\rm ML})^{\nu}}{\nu !}
  \left. \left[\frac{\partial^{\nu}}{\partial\xi^{\nu}}H(\xi^{\rm ML}|\xi )
    \right]\right|_{\xi =\xi^{\rm ML}}\, +\, R\, .
  \label{6.1}
\end{equation}
The quantity $R$ is the remainder. The zero-th order, $\nu =0\, ,$ of this expansion vanishes.
We expand up to $n=2$ and choose
Lagrange's remainder among several versions suggested in
Sect. 0.317 of \cite{Gradshteyn:15}.
This is   
\begin{equation}
  R\, =\,\frac{(\xi -\xi^{\rm ML})^3}{3!}
  \left.
  \frac{\partial^3}{\partial\xi^{\prime 3}}H(\xi^{\rm ML}|\xi^{\prime})
  \right|_{\xi^{\prime}=\xi^{\rm ML}+(\xi -\xi^{\rm ML})\Theta}\, ,
  \quad 0<\Theta <1\, .
  \label{6.2}
\end{equation}
The value of $\Theta$ remains open except for the fact that it lies between
$0$ and $1\, .$

\begin{figure}[htbp]
  \centering
  \includegraphics[width=12cm]{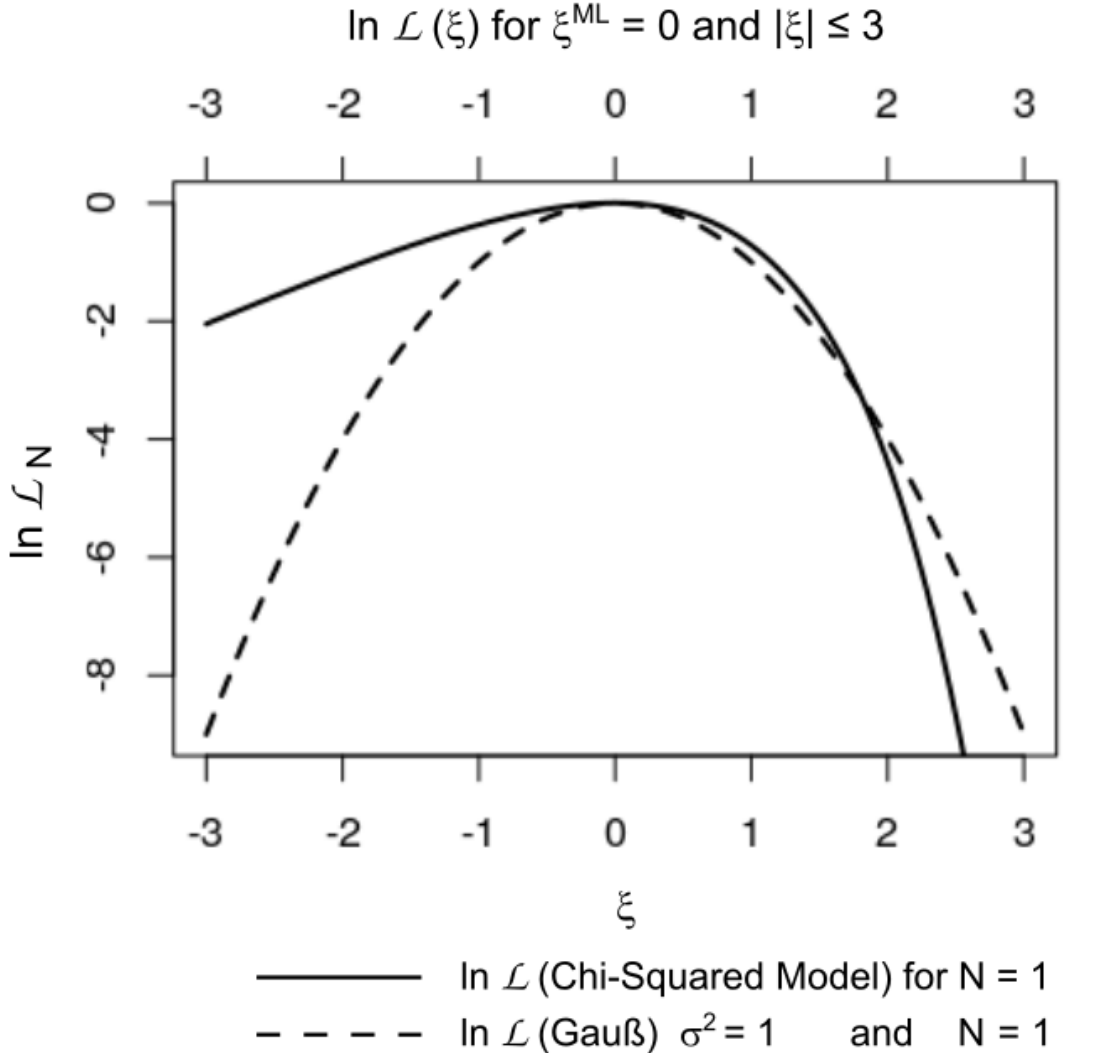}
  \caption{The graph compares a logarithmic likelihood function of a Chi-squared model
    to that of a Gaussian (i.e. a parabola). The former is taken from 
    Eq. (\ref{3.2}) together with (\ref{6.12c}) and $N=1\, ;$ the latter taken from (\ref{1.6f}) with
    $\sigma^2=F=1\, .$ The likelihood function of a Chi-squared model
    has a considerable skewness as compared to the Gaussian which is symmetric with respect to
    $\xi^{\rm ML}-\xi=0\, .$ The figure illustrates that skewness.
           } 
 \label{Gauss1} 
\end{figure}

After $N$ observations the Gaussian (\ref{1.6}) is accepted as a valid
approximation if the remainder $R$ is negligible for all $\xi -\xi^{\rm ML}$ in the interval
\begin{equation}
  -3\sigma \frac{1}{\sqrt{N}}\, <\, \xi -\xi^{\rm ML}\, <\, 3\sigma\frac{1}{\sqrt{N}}\, .
  \label{6.3}
\end{equation}
containing $99.73$ percent of the area of the Gaussian function.
We call it the $3\sigma$-interval of the Gaussian approximation (\ref{1.6}),
See FIG. \ref{Gauss1}. The value of $\sigma$ is given by the Fisher information 
\begin{equation}
  F\, =\, \sigma^{-2}
  \label{28a}
\end{equation}
of the model $p\, .$ This value is the expected curvature of the likelihood function.
The Fisher information together with the measure $\mu$ has been defined in Sect. \ref{sec:2.3}.

The term of $\nu =0$ in the expansion (\ref{6.1}) vanishes according to the definition (\ref{3.4})
of the functional $H\, .$ 
The first derivative $\frac{\partial}{\partial\xi}H$ vanishes at $\xi =\xi^{\rm ML}$
because $H$ attains its maximum value for $\xi =\xi^{\rm ML}\, .$
The second derivative $\frac{\partial^2}{\partial\xi^2}H$ at
$\xi =\xi^{\rm ML}$ yields the negative Fisher information according to
Eq. (\ref{4.3}). The third derivative --- in the remainder --- is calculated by help of the
translational invariance (\ref{3.4b}) of H. This gives
\begin{eqnarray}
  \left.\frac{\partial^3}{\partial\xi^{\prime 3}}H(\xi^{\rm ML}|\xi^{\prime})
  \right.|_{\xi^{\prime}=\xi^{\rm ML}+(\xi -\xi^{\rm ML})\Theta}
  \, &=&\,
  \left.\int_{-\infty}^{\infty}{\rm d}x\, p(x-\xi^{\rm ML})\frac{\partial^3}{\partial\xi^{\prime 3}}
  \ln p(x-\xi^{\prime})
  \right.|_{\dots} \nonumber\\
  \, &=&\,
  -\left.\int_{-\infty}^{\infty}{\rm d}x\, \frac{\partial^3}{\partial\xi^{\prime 3}}
  p(x-\xi^{\rm ML}+\xi^{\prime})\ln p(x)\right.|_{\dots} 
  \label{6.6}
\end{eqnarray}
The subscripts $|_{\dots}$ in the r.h.s. of this
equation mean repetitions of the subscript on the l.h.s.

In order to decide whether the remainder $R$ is small enough to be neglected, we must find below
the largest absolute value $|H^{(3)}|_{\rm max}$ of (\ref{6.6}). The largest value of $(\xi -\xi^{\rm ML})^3$
occurs at the upper end of the interval (\ref{6.3}), whence the maximum absolute value of the remainder is
\begin{equation}
  |R|_{\rm max}\, =\, \frac{1}{3!}\left(\frac{3\sigma}{\sqrt{N}}\right)^3\, |H^{(3)}|_{\rm max}\, .
  \label{6.6a}
\end{equation}
This shall be small compared to the second-order-term 
$$\frac{1}{2}\left(\frac{3\sigma}{\sqrt{N}}\right)^2\, F$$
everywhere in the interval (\ref{6.3}). Thus we require
\begin{equation}
  \frac{1}{2}\left(\frac{3\sigma}{\sqrt{N}}\right)^2\, F\, \gg\, \frac{1}{6}\left(\frac{3\sigma}{\sqrt{N}}\right)^3\, |H^{(3)}|_{\rm max}\, .
  \label{6.8}
\end{equation}
This gives, by use of (\ref{28a}),
\begin{equation}
  1\, \gg\, \frac{|H^{(3)}|_{\rm max}}{N^{1/2}F^{3/2}}\, .
  \label{6.9}
\end{equation}
We consider this condition to be fulfilled if
\begin{equation}
  0.1\, \ge\, \frac{|H^{(3)}|_{\rm max}}{N^{1/2}F^{3/2}}\, .
  \label{6.10}
\end{equation}
The last unequality is our criterion for the validity of the Gaussian approximation.
Replacing the requirement (\ref{6.9}) by (\ref{6.10}), we follow a convention. A step of one
``order of magnitude'' is usually considered to realise the requirement that one object be large
compared to another one.

\subsection{The Chi-Squared Model}
\label{sec:3.2}

As an example let us consider the model
\begin{equation}
  p(x-\xi )\, =\, \exp\left(x-\xi-e^{x-\xi}\right)
  \, ,\quad\quad -\infty <x,\xi <\infty\, .
  \label{6.11}
\end{equation}
It is normalised to unity, see section 8.312, no. 10 of \cite{Gradshteyn:15}.
This is a chi-squared distribution with two degrees of freedom transformed
such that it depends on the difference $x-\xi$ between event and
parameter, see Eq. (\ref{D.5}) in appendix \ref{D}. The ML estimator is $\xi^{\rm ML}=x\, ,$ see Eq. (\ref{B.4}) in appendix \ref{B}.

In the expansion (\ref{6.1}) the terms of order $\nu =0,1$ vanish according to Sect. \ref{sec:3.1}.
The derivatives of $H$ are obtained as follows. The first line of (\ref{3.4b}) with (\ref{6.11}) gives
\begin{eqnarray}
  \frac{\partial}{\partial\xi}H(\xi^{\rm ML}|\xi )\, &=&\, \frac{\partial}{\partial\xi}\int_{-\infty}^{\infty}{\rm d}x\, \Bigl(\exp (x-\xi^{\rm ML}-e^{x-\xi^{\rm ML}})\Bigr)
  \left[\xi^{\rm ML}-\xi-e^{x-\xi}+e^{x-\xi^{\rm ML}}\right]\nonumber\\
  &=&\, \int_{-\infty}^{\infty}{\rm d}x\, \Bigl(\exp (x-\xi^{\rm ML}-e^{x-\xi^{\rm ML}})\Bigr)\Bigl[-1+e^{x-\xi}\Bigr]\nonumber\\
  &=&\, -\int_{-\infty}^{\infty}{\rm d}x\, \exp\bigl(x-\xi^{\rm ML}-e^{x-\xi^{\rm ML}}\bigr)\nonumber\\
  & &\quad\quad\quad\, +\, \int_{-\infty}^{\infty}{\rm d}x\, \exp\bigl(x-\xi^{\rm ML}-e^{x-\xi^{\rm ML}}+x-\xi\bigr)\, . 
  \label{6.12}
\end{eqnarray}
In the last line the first term on the r.h.s. has the value $-1$ due to the normalisation of the model \ref{6.11}.
Whence, we obtain
\begin{equation}
  \frac{\partial}{\partial\xi}H(\xi^{\rm ML}|\xi )\, =\, -1\, +\, \int_{-\infty}^{\infty}{\rm d}x\, \bigl(\exp (x-\xi^{\rm ML}-e^{x-\xi^{\rm ML}}+x-\xi)\bigr)\, .
  \label{6.12aa}
\end{equation}
The substitution
\begin{equation}
  x^{\prime}\, =\, x-\xi^{\rm ML}
  \label{6.12a}
\end{equation}
yields
\begin{equation}
  \frac{\partial}{\partial\xi}H(\xi^{\rm ML}|\xi )\, =\, -1\, +\, \int_{-\infty}^{\infty}\, {\rm d}x^{\prime}\, \bigl(\exp (2x^{\prime}-e^{x^{\prime}}-\xi +\xi^{\rm ML})\bigr)\, .
  \label{6.12c}
\end{equation}

For $\xi =\xi^{\rm ML}$ this becomes
\begin{eqnarray}
 \frac{\partial}{\partial\xi}H(\xi^{\rm ML}|\xi )\left|_{\xi=\xi^{\rm ML}}\right.\, &=&\, -1\, +\, \int_{-\infty}^{\infty}\, {\rm d}x^{\prime}\, \exp\bigl(2x^{\prime}-e^{x^{\prime}}\bigr)\nonumber\\
  &=&\, -1\, +\, \Gamma (2)\, ,
  \label{6.12cc}
\end{eqnarray}
see appendix \ref{E}. This yields
\begin{eqnarray}
  \frac{\partial}{\partial\xi}H(\xi^{\rm ML}|\xi )\left|_{\xi =\xi^{\rm ML}}\right.\, &=&\, -1\, +\, 1\nonumber\\
  &=&\, 0\, .
  \label{6.12cd}
\end{eqnarray}
This result was expected since $H(\xi^{\rm ML}|\xi )$ has an extreme value at $\xi =\xi^{\rm ML}\, .$

The second derivative of $H$ with respect to $\xi$ is
\begin{eqnarray}
  \frac{\partial^2}{\partial\xi^2}H(\xi^{\rm ML}|\xi )\, &=&\, -\int_{-\infty}^{\infty}\, {\rm d}x^{\prime}\, \exp\bigl(2x^{\prime}-e^{x^{\prime}})\bigr)e^{-\xi +\xi^{\rm ML}}\nonumber\\
  &=&\, -e^{\xi^{\rm ML}-\xi}\, \int_{-\infty}^{\infty}\, {\rm d}x^{\prime}\, \exp\bigl(2x^{\prime}-e^{x^{\prime}})\bigr)\nonumber\\
  &=&\, -e^{\xi^{\rm ML}-\xi}\, \Gamma (2)\nonumber\\
  &=&\, -e^{\xi^{\rm ML}-\xi}\, .
  \label{6.12ce}
\end{eqnarray}
The definition of $H(\xi^{\rm ML}|\xi )$ in the first line of Eq. (\ref{3.4b}) together with the substitution (\ref{6.12a}) yields
a symmetry relation of the functional $H\, ,$
\begin{eqnarray}
  H(\xi^{\rm ML}|\xi )\, &=&\, \int_{-\infty}^{\infty}\, {\rm d}x^{\prime}\, \exp\bigl(x^{\prime}-e^{x^{\prime}})\bigr)
  \ln\frac{\exp\bigl(x^{\prime}-e^{x^{\prime}+\xi^{\rm ML}-\xi}+\xi^{\rm ML}-\xi\bigr)}{\exp\bigl(x^{\prime}-e^{x^{\prime}}\bigr)}\nonumber\\
    &=&\, H(0|\xi -\xi^{\rm ML})\, .
    \label{6.12cf}
\end{eqnarray}

The last line of Eq. (\ref{6.12ce}) yields the derivatives
\begin{equation}
  \frac{\partial^{\nu}}{\partial\xi^{\nu}}H(\xi^{\rm ML}|\xi )\, =\, (-1)^{\nu +1}e^{\xi^{\rm ML}-\xi}\, ,\quad\quad\quad \nu\ge 2\, .
      \label{6.14}
\end{equation}
One especially obtains the Fisher information
\begin{eqnarray}
  F\, &=&\, -\left.\frac{\partial^2}{\partial\xi^2}H(0|\xi -\xi^{\rm ML})\right|_{\xi =\xi^{\rm ML}}\nonumber\\
      &=&\, 1\, .
  \label{6.15}
\end{eqnarray}

Equation (\ref{6.14}) gives the maximum $|H^{(3)}|_{\rm max}$ of the absolute value of
(\ref{6.6}) in the interval (\ref{6.3}) to be
\begin{eqnarray}
  |H^{(3)}|_{\rm max}\,
  &=&\, \exp\left(3\frac{\sigma}{N^{1/2}}\right)\nonumber\\
  &=&\, \exp\left(\frac{3}{F^{1/2}N^{1/2}}\right)\nonumber\\
  &=&\, \exp\left(\frac{3}{N^{1/2}}\right)\, .
  \label{6.16}
\end{eqnarray}
Then Eq. (\ref{6.10}) yields the criterion
\begin{equation}
  0.1\, \ge\, \frac{1}{N^{1/2}}\exp \left(\frac{3}{N^{1/2}}\right)\, .
      \label{6.17}
\end{equation}

Table \ref{tab:1} shows that the unequality (\ref{6.17}) is satisfied if
\begin{equation}
  N\, \ge\, 160\, .
  \label{6.20}
\end{equation}
Thus the posterior of the chi-squared distribution (\ref{6.11}) with two
degrees of freedom can be considered Gaussian when the number $N$ of events is
larger than or equal to $160\, .$
This shows that an intuition based on the central limit theorem, were a current ``rule of thumb'' says
$N\approx 30\, ,$ can be misleading. See Sect. 7.4.2 of \cite{Ross:17} or Hogg et al. \cite{Hogg:15}.

\begin{table}[t]
  \begin{tabular}{|r|r|}
  \multicolumn{2}{l}{Table I}\\
\hline
$N$    &  r.h.s. of (\ref{6.17})\\
\hline
    3   &  3.263\\
    4   &  2.241\\
    5   &  1.711\\
   10   &  0.817\\
   20   &  0.437\\
   30   &  0.316\\
   40   &  0.254\\
   50   &  0.216\\
   75   &  0.163\\
  100   &  0.135\\
  150   &  0.104\\
  155   &  0.102\\
  160   &  0.100\\
  165   &  0.098\\
  \hline
\end{tabular}
  \caption{Right hand side of inequality (\ref{6.17}) for various values of $N$.}
\label{tab:1}
\end{table}

To find the necessessary $N\, ,$ we have not used the article by Berry \cite{Berry:41}
on the accuracy of the {G}aussian approximation because Berry does not discuss the measure on the scales of the variables.

\section{The Binomial Model}
\label{sec:4}

In the following Sects. \ref{sec:4} and \ref{sec:5} the
{\it binomial model} $q(x|\xi )$ is considered. It is also called the
model of the ``simple alternative'' since its event $x$ is restricted to two
possible values $x=0$ or $x=1$ while $\xi$ is a continuous variable which conditions
the probability to find $x\, .$ By consequence, form
invariance cannot be expressed by the difference between $x$ and $\xi\, .$
Yet form invariance under translations is found by reformulating the model
such that it depends on the difference between the ML
estimator $\xi^{\rm ML}$ and the parameter $\xi\, ,$ and only on this difference.

\subsection{Definition of the Binomial Model}
\label{sec:4.1}

The binomial model is given by
\begin{equation}
  q(x|\xi )\, =\, \left[R(\xi )\right]^x\left[1-R(\xi )\right]^{1-x}\, ,
  \quad x=0,1\, .
  \label{7.1}
\end{equation}
Here, $x=1$ can be interpreted as the ``correct'' and $x=0$ as the ``false''
answer to a given question. Instead of correct or false, the two possibilities
may also mean one or the other side of a thrown coin. The 
function $R$ of the continuous variable $\xi$ is called item response function.

The name of {\it item response function} (IRF) is due to authors who discussed
the ideas of  ``item response theory'' \cite{Thurstone:25,Ferguson:42,Lazarsfeld:73,Lord:80,Rasch:60,Rasch:67,Rasch:77,Rasch:80}
applied in intelligence tests as well as educational tests, e.g. the PISA tests \cite{PISA:15}.
Differing from them we define the IRF by the requirement that the measure $\mu $ on the scale
of $\xi$ should be constant. This means that the Fisher information $F$ shall be independent of $\xi\, ,$
see Eq. (\ref{4.4a}). The Fisher information of the binomial model is defined much as in Eq. (\ref{1.3}),
except for the integral over $x$ which becomes a sum over $x=1,2\, .$
We make use of the second version of Eq. (\ref{4.1}) and obtain a differential equation with the solution
\begin{equation}
  R(\xi )\, =\, \cos^2\xi\, ,\quad\quad -\pi /2\le \xi\le \pi /2\, ,
  \label{7.3}
\end{equation}
see appendix \ref{J}.
The domain of definition of the parameter $\xi$ makes sure that the likelihood function possesses
a maximum and vanishes towards the ends of the domain. This is analogous to the model $p(x-\xi )$
introduced in Sect. \ref{sec:2.1}. The binomial model with the IRF (\ref{7.3}) is called the
{\it trigonometric model}.

By the definition of the IRF (\ref{7.3}) the Fisher information of the model (\ref{7.1}) becomes
\begin{equation}
  F\, \equiv\, 4\, ,
  \label{9.4}
\end{equation}
see appendix \ref{J};
whence the measure (\ref{4.4a}) on the scale of $\xi$ is
\begin{equation}
  \mu\, \equiv\, 2\, .
  \label{9.1a}
\end{equation}

G. Rasch did not consider the measure on the scale of the IRF, he derived
the IRF by requiring the property of ``specific objectivity'' to any measurement
\cite{Rasch:60,Rasch:67,Rasch:77,Rasch:80}. By this property the score of correct answers
became the sufficient statistic in his model. The trigonometric IRF, however, is compatible with specific objectivity,
see Sect. 3.1 of \cite{Fuhrm:16}.

For $N$ events $\bold{x}=(x_1,\dots ,x_N)$ conditioned by one and the same
$\xi$, the binomial model reads
\begin{equation}
  q_N(\bold{x}|\xi )\, =\, \prod_{k=1}^N\,
  \left[R(\xi )\right]^{x_k}\left[1-R(\xi )\right]^{1-x_k}\, ,
  \quad x_k=0,1\, .
  \label{7.4}
\end{equation}
This can be rewritten by help of the score $s_c$ of the answers that yield $x=1$ as well as the number $N-s_c$
of the answers that yield $x=0\, .$ One obtains
\begin{equation}
  q_N(s_c|\xi )\, =\, \Big(
  \begin{array}{c}
    N\\
    s_c
  \end{array}
  \Big)\left[R(\xi )\right]^{s_c}\left[1-R(\xi )\right]^{N-s_c}\, ,
  \quad 0\le s_c\le N\, ,\, N\ge 1\, .
  \label{7.5}
\end{equation}
The quantity
$\Big(
  \begin{array}{c}
    N\\
  s_c
  \end{array}
  \Big)$
is a binomial coefficient. By applying the binomial formula one finds that
$q_N$ is normalised,
\begin{equation}
    \sum_{s_c =0}^N\, q_N(s_c|\xi )\, =\, 1\, .
    \label{7.7}
\end{equation}
The posterior distribution will be obtained in Sect. \ref{sec:4.5}.

\subsection{The ML Estimator for the Trigonometric Model}
\label{sec:4.2}

Given the event $s_c$ in the framework of the distribution (\ref{7.5}), the ${\rm ML}$ estimator
is found by solving the the ML equation
\begin{equation}
  \frac{\partial}{\partial\xi}\ln q_N(s_c|\xi )\, =\, 0\, .
  \label{8.1}
\end{equation}
This leads to
\begin{eqnarray}
  0\, &=&\, \frac{\partial}{\partial\xi}\left[s_cR(\xi )+(N-s_c)\ln\big(1-R(\xi )\big)\right]\nonumber\\
  &=&\, s_c\frac{R^{\prime}}{R}\, -\, (N-s_c)\frac{R^{\prime}}{1-R}\, .
  \label{8.2}
\end{eqnarray}
With the IRF (\ref{7.3}) this ML-equation becomes
\begin{eqnarray}
  0\, &=&\, -2\cos\xi\sin\xi\left[\frac{s_c}{\cos^2\xi}-\frac{N-s_c}{\sin^2\xi}\right]\nonumber\\
  &=&\, -2\frac{s_c-N\cos^2\xi}{\cos\xi\sin\xi}
  \label{8.4}
\end{eqnarray}
which is solved by $\xi^{\rm ML}(s_c)$ such that
\begin{equation}
  \cos^2\xi^{\rm ML}\, =\, \frac{s_c}{N}\, .
  \label{8.5}
\end{equation}
The denominator of the r.h.s. of (\ref{8.4})  seems to exclude the values $\xi^{\rm ML}=\pm\pi /2$ and $\xi^{\rm ML}=0$
because it vanishes there. These values correspond to the ``uniform answers'' $s_c=0$ and $s_c=N\, ,$ respectively.
We show that the uniform answers are not excluded.
In the case of $s_c=0$ equation (\ref{8.4}) turns into
\begin{equation}
  0\, =\, \frac{\cos\xi}{\sin\xi}
  \label{8.6}
\end{equation}
and has the solutions $\xi^{\rm ML}=\pm\pi /2\, .$ This is conform with (\ref{8.5}).
In the case of $s_c=N$ equation (\ref{8.4}) becomes
\begin{eqnarray}
  0\, &=&\, \frac{1-\cos^2\xi}{\cos\xi\sin\xi}\nonumber\\
  &=&\, \frac{\sin\xi}{\cos\xi}
  \label{8.7}
\end{eqnarray}
which is solved by $\xi^{\rm ML}=0$ and again conforms with Eq. (\ref{8.5}).

\subsection{The Likelihood Function of the Trigonometric Model in the case of $N=1$}
\label{sec:4.3}

The likelihood function of the trigonometric model $q_{N=1}$ with $R$ according to (\ref{7.3})
is form invariant. To show this we study again the two cases of $s_c=0,1\, .$

(i) Let $s_c=0$ be observed. Then the likelihood function is proportional to $1-R(\xi )$
and the ML estimator is $\xi^{\rm ML}=\pm\pi /2\, .$ Thus the likelihood function is
\begin{equation}
  {\cal L}_{N=1}(\xi )\, \propto\, \sin^2\xi\, ,\quad\quad\quad -\pi /2<\xi <\pi /2\, .
  \label{10.4}
\end{equation}
This can also be expressed as
\begin{equation}
  {\cal L}_{N=1}(\xi )\, \propto\, \cos^2(\xi -\xi^{\rm ML})\, .  
  \label{10.4a}
\end{equation}
(ii) Now let $s=1$ be observed. Then the likelihood function is proportional to $R(\xi )\, $ whence
the ML estimator is $\xi^{\rm ML}=0\, ,$ according to the case of Eq. (\ref{8.7}).
Thus the likelihood function is
\begin{equation}
  {\cal L}_{N=1}(\xi )\, \propto\, \cos^2\xi\, .
  \label{8.3}
\end{equation}
Since now $\xi^{\rm ML}=0$ this is again expressed by Eq. (\ref{10.4a}).

Thus if $N=1$ and $R$ given by (\ref{7.3}) the likelihood function
of the trigonometric model depends on the difference $\xi -\xi^{\rm ML}$ and only on this difference.
This is what we call form invariance under translations. In summary: The posterior of the trigonometric model $q$ is
obtained with a constant prior and it
is form invariant for $N=1\, .$

The result (\ref{10.4a}) is found for every value of $\xi^{\rm ML}$ given by Eq. (\ref{8.5}).
Thus it remains true for any number $N$ of observations. The posterior of the trigonometric model is form invariant
under translations in the sense that it depends on the difference $\xi -\xi^{\rm ML}\, .$ However, the variables
$\xi$ and $\xi^{\rm ML}$ do not ``make the same use'' of the interval $-\pi /2<\xi <\pi /2$ given in Eq. (\ref{10.4}):
The variable $\xi$ is defined everywhere on this interval. The ML estimator $\xi^{\rm ML}$ assumes only a finite
number of values within that interval. In the following Sect. \ref{sec:4.4} the trigonometric model itself
--- not only its posterior --- will be brought into translational form invariance.

\subsection{The Trigonometric Model with Translational Invariance}
\label{sec:4.4}

For arbitrary $N$ the ML estimator given by Eq. (\ref{8.5}) is not restricted to two values.
For sufficiently large $N$ any value in the interval $-\pi /2\le\, \xi^{\rm ML}\, \le 0$ can be
approached arbitrarily closely. Is it possible to define the binomial model such that the model
itself is translationally form invariant? Yes, this is possible and leads to the same Fisher information (\ref{9.4}).

Consider the model
\begin{equation}
  t(s|\xi )\, =\, {\cal N}\cos^2(s-\xi )\, ,\quad\quad
  -\pi /2\, \le\, s,\xi\, \le\pi/2\, ,
  \label{10.10}
\end{equation}
which depends on the difference between the observation $s$ and the condition $\xi\, .$

According to the argument Sect. \ref{4.3} on the posterior of the trignometric model,
the ML estimator comes arbitrarily close to every value in the interval
$-\pi /2\le \xi^{\rm ML}\le \pi /2$ since, for arbitrary $N\, ,$ it assumes every rational number
in this interval. Thus the trigonometric model can be extended to the model \ref{10.10}.
The fact that $\xi^{\rm ML}$ is defined within, $-\pi /2\le \xi^{\rm ML}\le \pi /2$ means that
$\xi$ is defined in the same interval.

We shall show that (\ref{10.10}) is translationally form invariant and has the Fisher information
(\ref{9.4}) in agreement with the trigonometric model (\ref{7.1}) with (\ref{7.3}).

The normalisation ${\cal N}$ in (\ref{10.10}) equals the integral over the domain of $s\, ;$
it is independent of $\xi$ because the domain given in (\ref{10.10}) is an interval of length $\pi$ which is equal
to one period of the periodic function $\cos^2(x-\xi)\, .$ Thus ${\cal N}$ is given by the integral
\begin{eqnarray}
  {\cal N}^{-1}\, &=&\, \int_{-\pi /2}^{\pi /2}{\rm d}s\, \cos^2(s-\xi )\nonumber\\
  &=&\, \int_{-\pi /2}^{\pi /2}{\rm d}s\, \cos^2s\, .
  \label{10.14}
\end{eqnarray}
One obtains
\begin{equation}
  {\cal N}\, =\, 2/\pi\, ,
  \label{10.14b}
\end{equation}
see appendix \ref{G}.

The ML estimator $\xi^{\rm ML}$ for a given event $s$ occurs at
\begin{equation}
  \xi^{\rm ML}\, =\, s
  \label{10.26d}
\end{equation}
since the likelihood function $t(s|\xi )$ of (\ref{10.10}) becomes maximal when the argument
of the $\cos^2$-function is zero.

The posterior of $t(s|\xi )$ shall be a function of $\xi^{\rm ML}-\xi$ as it is in the context
of the model (\ref{6.11}) in Sect. \ref{sec:3}. This requires to shift --- within the posterior ---
the value of $\xi^{\rm ML}$ to the value of zero. We can do so since the measure on the scale of
$\xi$ is constant. The shift avoids that the distance $|\xi^{\rm ML}-\xi|$ bacomes larger than the
length $\pi$ of the interval in which $\xi$ is defined. In the context of the model (\ref{6.11})
such a shift was not needed since every difference $\xi^{\rm ML}-\xi$ was contained in the infinite
domain $-\infty <\xi < \infty\, ,$ where $\xi$ was defined.

Since $\xi^{\rm ML}$ is the sufficient statistic of the
model (\ref{10.10}) the score $s_c$ of correct answers determines the sufficient statistic
when $N$ is given. In this sense the model (\ref{10.10}) confirms the requirement of
G. Rasch \cite{Rasch:60,Rasch:66a,Rost:04,Fisch:81,Fisch:95}
that the score should be the sufficient statistic of the binomial model. 
See especially Chap. 4.6 of \cite{Fuhrm:16} and Chap. 12.3.1 of \cite{Harney:16}, were so-called
Guttman schemes \cite{Guttman:73} are analysed.

The Fisher information of the model (\ref{10.10}) is
\begin{equation}
  F(\xi )\, =\, -{\cal N}\int_{-\pi /2}^{\pi /2}{\rm d}s\,
  \cos^2(s-\xi )\frac{\partial^2}{\partial\xi^2}\ln\cos^2(s-\xi)
  \label{10.17}
\end{equation}
by Eq. (\ref{4.1}). For the second derivative in the integrand
one finds
\begin{equation}
  \frac{\partial^2}{\partial\xi^2}\ln\cos^2(s-\xi)\, =\, -\frac{2}{\cos^2(s-\xi )}\, .
  \label{10.17a}
\end{equation}
Together with (\ref{10.14b}) this gives
\begin{equation}
  F(\xi )\, \equiv\, 4\, .
  \label{10.26}
\end{equation}
Thus the Fisher information of the translationally invariant model (\ref{10.10}) agrees with the Fisher information
(\ref{9.4}) obtained from the binomial model (\ref{7.1}).

\subsection{The Functional $H$ for the Trigonometric Model with Translational Invariance}
\label{sec:4.5}

We express the logarithmic likelihood function $\ln{\cal L}_N$ for $N$ events
given by the model (\ref{10.10}) in analogy to Eq. (\ref{3.2}) as
\begin{eqnarray}
  \ln{\cal L}_N\, &=&\, \, N\, H(\xi^{\rm ML}|\xi )\nonumber\\
                  &=&\, N\, \int_{-\pi /2}^{\pi /2}{\rm d}s\, t(s-\xi^{\rm ML})\ln\frac{t(s-\xi )}{t(s-\xi^{\rm ML})}\, .
  \label{10.11}
\end{eqnarray}
The integration extends over one period
of the $\cos^2$-function. Inserting (\ref{10.10}) into (\ref{10.11})
one obtains
\begin{equation}
  \ln{\cal L}_N\, =\, N\, {\cal N} \int_{-\pi /2}^{\pi /2}{\rm d}s\,
                   \cos^2(s-\xi^{\rm ML})\ln\frac{\cos^2(s-\xi )}{\cos^2(s-\xi^{\rm ML})}\, .
  \label{10.12a}
\end{equation}

\begin{figure}[htbp]
  \centering
  \includegraphics[width=12cm]{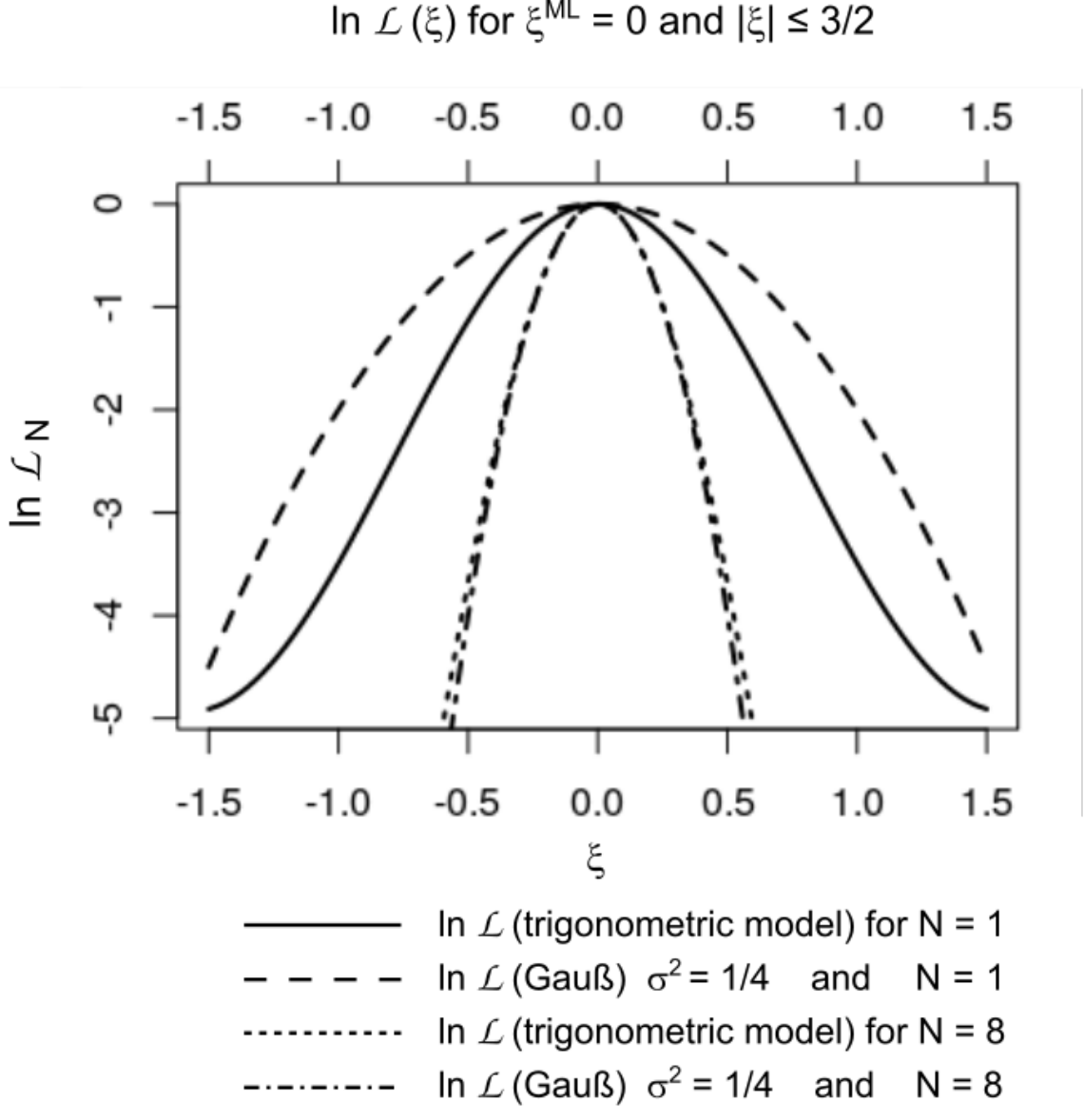}
  \caption{The graph compares a logarithmic likelihood function of the trigonometric model with form invariance to that of a Gaussian
    (i.e. a parabola) for $N=1\, $ and $N=8\, $. The former is taken from Eq. (\ref{10.12a}), the latter from
    Eq. (\ref{1.6f}) with $\sigma^2=F^{-1}= 1\, .$
    The curvature of the Gaussian agrees with the curvature of the likelihood function at their maxima.
    This figure also illustrates the reflection symmetry of the trigonometric model with respect to $\xi =0\, .$
    With increasing $N$ the likelihood functions are ever more concentrated around their maximum value.
    This is illustrated by the comparison between $N=1$ and $N=8\, .$
    } 
 \label{Gauss2} 
\end{figure}

The functional $H(\xi^{\rm ML}|\xi )$ exists and has the property (\ref{3.4b}) because the integral in (\ref{10.12a}) exists
although the function $\cos^2s^{\prime}$ with, e.g.,
\begin{equation}
  s^{\prime}\, =\, s-\xi
  \label{10.12b}
\end{equation}
will vanish at an isolated point $s_0^{\prime}$ and $\ln\cos^2s_0^{\prime}$ will diverge. Appendix \ref{H} shows that the integral
over an $\epsilon$-interval that includes $s_0^{\prime}$, does exist. Here, $\epsilon$ may be arbitrarily small. The value of
this integral is not proportional to $\epsilon\, ;$
instead it is proportional to $\epsilon\ln\epsilon\, .$ Although
this is not as small as $\epsilon\, ,$ it approaches zero when $\epsilon\to 0\, .$
Thus the logarithmic likelihood $\ln{\cal L}_N$ of Eq. (\ref{10.12a}) exists
even if the integration runs over a point where $\ln\cos^2$ diverges.

Applying the substitution (\ref{10.12b}) to the integrand in (\ref{10.12a}) one obtains
\begin{eqnarray}
  \ln{\cal L}_N\, &=&\, N{\cal N}\int_{-\pi /2}^{\pi /2}{\rm d}s^{\prime}\, \cos^2s^{\prime}\ln\frac{\cos^2(s^{\prime}-\xi +\xi^{\rm ML})}{\cos^2s^{\prime}}\nonumber\\
    &=&\, NH(0|\xi -\xi^{\rm ML})\, .
  \label{10.12k}
\end{eqnarray}
The substitution does not require a shift of the limits of integration because the
integrand is periodic with a period of length $\pi$ and the integration
covers one period. Thus the functional $H(\xi^{\rm ML}|\xi )$ for the translationally invariant
trigonometric model has the property (\ref{3.4b}) which was found earlier for
models defined on the entire real axis.

The fact that ${\cal L}_N$ equals $N$ times $H(0|-\alpha )$ shows that Eq. (\ref{5.4}) holds here, too: The posterior distribution from $N$
events of the trigonometric model is proportional to the $N$-th power of the posterior from one event.

\section{The {G}aussian Approximation to the Trigonometric Model with Translational Invariance}
\label{sec:5}
 
The present section establishes the criterion for the validity of the Gaussian
approximation to the likelihood function of the trigonometric model with translational invariance.
In analogy to the procedure in Sect. \ref{sec:3} we expand $H(\xi^{\rm ML}|\xi )$ into a Taylor series
with respect to $\xi$ at $\xi^{\rm ML}\, .$ Again the zero-th order of this expansion vanishes by the
definition of $H$ in Eq. (\ref{3.4}). Again the likelihood function is related to $H$ via Eq. (\ref{3.2}).
Using the abbreviation (\ref{F.4}), introduced in appendix \ref{F}, we obtain
\begin{equation}
  H(0|\xi -\xi^{\rm ML})\, =\, H(\xi^{\rm ML}-\xi |0)\, .
  \label{11.0}
\end{equation}

Appendix \ref{F} furthermore shows that the expression (\ref{11.0}) is mirror symmetric with respect to $\xi^{\rm ML}-\xi =0\, ,$ i.e.
\begin{equation}
  H(0|\xi -\xi^{\rm ML})\, =\, H(0|\xi^{\rm ML}-\xi )\, .
  \label{11.0a}
\end{equation}
Therefore all odd derivatives of $H(0|\xi^{\rm ML}-\xi )$ vanish at $\xi^{\rm ML}-\xi =0\, ,$  thus
\begin{equation}
  \left.\frac{\partial^{2\nu +1}}{\partial\xi^{2\nu +1}}\, H(0|\\xi^{\rm ML}-\xi )\right|_{\xi =\xi^{\rm ML}}\, =\, 0\, ,
  \quad\quad\quad \nu =\, 0,1,2,\dots\, .
  \label{11.1}
\end{equation}
By consequence, the remainder of the present expansion is not of the
third order as it is in Sect. \ref{sec:3.1}; it rather is of the fourth order.

All even derivatives $\frac{\partial^{2\nu}}{\partial\alpha^{2\nu}}H(0|\alpha )$ exist. Therefore in the present context the Taylor expansion (\ref{6.1}) becomes 
\begin{equation}
  H(\xi^{\rm ML}|\xi )\, =\, -\frac{(\xi^{\rm ML}-\xi )^2}{2}F\, +\, \frac{(\xi^{\rm ML}-\xi )^4}{4!}\left.\frac{\partial^4}{\partial\xi^4}H(\xi^{\rm ML}|\xi )\right|_{\xi =\xi^{\rm ML}+(\xi^{\rm ML}-\xi )\Theta}\, .
  \label{11.2}
\end{equation}

In analogy to Sect. \ref{sec:3.1} the Gaussian approximation is considered
valid when the term of fourth order is negligible as compared to the term of
second order for all $\xi -\xi^{\rm ML}$ in the interval (\ref{6.3}). See FIG. \ref{Gauss2}.
For these values of $\xi$ we have to find the maximum of the absolute value of the
fourth derivative in Eq. (\ref{11.2}). We call it $|H^{(4)}|_{\rm max}\, .$
In analogy to Eq. (\ref{6.8}) the Gaussian approximation is accepted when
\begin{equation}
\frac{1}{2}\big(\frac{3\sigma}{\sqrt{N}}\big)^2F\, \gg\, \frac{1}{24}\big(\frac{3\sigma}{\sqrt{N}}\big)^4|H^{(4)}|_{\rm max}
\label{11.3}
\end{equation}
or
\begin{equation}
  1\, \gg\, \frac{3}{4}\frac{\sigma^2}{N}\frac{|H^{(4)}|_{\rm max}}{F}\, .
  \label{11.3b}
\end{equation}
As in Eq. (\ref{28a}) the Fisher information $F$ equals $\sigma^{-2}$ and we obtain
the condition
\begin{equation}
  1\, \gg\, \frac{3|H^{(4)}|_{\rm max}}{4NF^2}\, .
  \label{11.3a}
\end{equation}
In appendix \ref{F} the value of $|H^{(4)}|_{\rm max}$ is found to be
\begin{equation}
  |H^{(4)}|_{\rm max}\, =\, 16\, .
  \label{11.4}
\end{equation}
The value of the Fisher information $F$ is given by Eq. (\ref{10.26}). Whence,
the condition (\ref{11.3a}) for the Gaussian approximation becomes
\begin{equation}
  1\, \gg\, \frac{3}{4N}\, .
  \label{11.5}
\end{equation}
We consider it to be fulfilled if
\begin{equation}
  0.1\, \ge\, \frac{3}{4N}
  \label{11.6}
\end{equation}
or
\begin{equation}
  N\, \ge\, 8\, .
  \label{11.7}
\end{equation}

Thus the condition (\ref{11.7}) for the Gaussian approximation to the
binomial model is much more easily fulfilled than the corresponding condition
(\ref{6.10}) for the chi-squared model. The reason is:
In Sect. \ref{sec:3.1} the remainder was of third order;
in the present case it is of fourth order.

The step from (\ref{11.5}) to (\ref{11.6}) is due to the common idea that a given
positive $x$ be large against $y>0$ when $x$ is one order of magnitude larger than
$y\, .$ This idea is commonly used for approximations in mathematics \cite{BSMM:07}
and physics \cite{CRLC:04}. If a higher accuracy is required, this rule can
easily be adapted.

In Chap. 12 of \cite{Harney:16} a version of Item Response Theory is presented
which makes use of the trigonometric IRF (\ref{7.3}). In that context
a (simulated) competence test was discussed which asked $20$ questions. From
the present result (\ref{11.7}) follows that the
estimated person-parameters have a Gaussian posterior distribution.

In Sect. \ref{sec:3.2} as well as in the present section the condition for the
Gaussian approximation is independent of the value of $\xi^{\rm ML}\, .$
In the case of the binomial model, this means that the Gaussian approximation
is valid even for uniform or close to uniform patterns of answers.

\section{Conclusions}
\label{sec:6}

The present text presents a criterion for the validity of the Gaussian
approximation to the likelihood function of a statistical model
$p$ when $N$ observations $\bold{x}=(x_1,\dots ,x_N)$ have been
collected.

We require that a statistical model possesses a symmetry between the observed quantity $x$
and the parameter $\xi\, .$ This symmetry is defined in terms of a Lie group. It has been
called form invariance. It allows to formally specify the prior distribution ---
required by Bayes but not specified by him. The model $p$ can then be parameterised such
that it depends on the difference between the observed quantity $x$ and the parameter $\xi$
which conditions the observations. The model $p(x-\xi )$ remains invariant when both quantities
are shifted by the same amount. We call this property translational form invariance. Then
the likelihood function ${\cal L}_N(\xi )$ is shown to depend on the difference $\xi^{\rm ML}-\xi $
between the maximum likelihood estimator $\xi^{\rm ML}$ and the parameter $\xi\, .$
This means that the Bayesian posterior distribution, too, depends on $\xi^{\rm ML}-\xi\, .$
One can even shift the scale of $\xi$ such that the ML estimator is found at $\xi^{\rm ML}=0\, .$

A total of $N$ observations conditioned by one and the same parameter $\xi$ is generally expected
to lead to a Gaussian likelihood function for sufficiently large $N\, .$ The question of how large
$N$ must be, in order to justify the Gaussian approximation, is answered in Sects. \ref{sec:3}
and \ref{sec:5} for two quite different examples. The basic idea is, that a valid Gaussian
approximation to the posterior distribution means that the error assigned to the parameter
$\xi$ equals $\sigma = (NF)^{-1/2}$ within the interval $|\xi^{\rm ML}-\xi |< 3\sigma\, .$
Here, $F$ is the Fisher information yielded by the model $p\, .$ The probability for
$\xi$ to lie outside this $3\sigma$-interval is neglected. A stricter condition for the
Gaussian approximation would be achieved by requiring an interval larger than $3\sigma$ for its validity. 

The example of Sect. \ref{sec:3} is a version of the chi-squared model with
two degrees of freeedom. This distribution, as well as its posterior,
has a considerable skewness. In this case one
needs $N=160$ observations for the Gaussian approximation to be acceptable. 
The large number of required observations is due to the skewness.

The example of Sect. \ref{sec:5} is the trigonometric model, a specific form of the binomial model based on form invariance. It strongly differs
from the chi-squared model since the observations are not taken from a continuum
of real numbers but rather from the alternative of $0$ or $1\, .$ The likelihood
function turns out to exhibit a mirror symmetry --- as the Gaussian exhibits, too. This 
helps to approach the Gaussian distribution. We find:
For $N\ge 8$ the posterior of the trigonometric model can be considered as {G}aussian.
This holds for every value of the ML estimator. It is a favorable result for the
application of the binomial model to competence tests such as the PISA studies
\cite{PISA:15}.

In general, we see a practical interest in our results since the normal distribution
is the basis of parametric methods in applied statistics, widely used in many areas
(education, medicine, science, etc.). To know whether the normal distribution is
applicable or not, is of interest for practitioners in these fields.

\begin{appendix}
\section{Comparing Two Distributions}
\label{A}

We show that for any two normalised distributions $p$ and $q$ the unequality 
\begin{equation}
  \sum_l\, p_l\ln \frac{q_l}{p_l}\, \le\, 0
  \label{A.1}
\end{equation}
holds provided that $p$ and $q$ are labelled by the entire numbers $l$ and are
normalised according to
\begin{eqnarray}
  \sum_l\, q_l\, &=&\, 1\, ,\nonumber\\
  \sum_l\, p_l\, &=&\, 1\, .
  \label{A.2}
\end{eqnarray}
The l.h.s. of (\ref{A.1}) attains its maximum value of zero when and only when $p_l$ and
$q_l$ agree for every $l\, .$

The unequality (\ref{A.1}) is a consequence of the unequality
\begin{equation}
  \ln s\, \le\, s-1\,  \quad\quad\quad {\rm for}\,\, \, s>0\, .
  \label{A.3}
\end{equation}
The linear function $s-1$ is tangent to the function $\ln s\, .$
The common point lies at $s=1\, ,$ where both functions have the value of $1\, .$
Setting $s=q_l/p_l$ the unequality entails
\begin{equation}
  \ln \frac{q_l}{p_l}\, \le\, \frac{q_l}{p_l}-1 
  \label{A.4}
\end{equation}
or
\begin{equation}
  p_l\ln \frac{q_l}{p_l}\, \le\, q_l-p_l 
  \label{A.5}
\end{equation}
for every $l\, .$ The quantity on the l.h.s. has also been introduced by Campbell
in Sect. 1 of Chap. 5 of \cite{Karmeshu:03}. Summing (\ref{A.5}) over $l$ yields the unequality (\ref{A.1}).
When the distributions $p$ and $q$ agree whith each other, the l.h.s. of
(\ref{A.1}) vanishes. Then and only then the expression assumes its maximum value.

One can interprete $p_l$ as the probability $p(x_l-\xi^{\rm ML})\frac{\Delta x}{N}$
contained in a bin centered at the value $x_l$ of the real variable $x$ and
having the width $\frac{\Delta x}{N}\, .$ Here, $p$ shall be a normalised probability
density. Similarly one can interprete $q_l$ as the probability
$p(x_l-\xi )\frac{\Delta x}{N}\, .$ In the limit of $\frac{\Delta x}{N}\to 0$ the unequality
(\ref{A.1}) then yields the unequality
\begin{equation}
  \frac{1}{N}\int_{{\rm domain}\, x}{\rm d}x\,
  p(x-\xi^{\rm ML})\ln\frac{p(x-\xi )}{p(x-\xi^{\rm ML})}\, \le\, 0 .
    \label{A.6}
\end{equation}
or
\begin{equation}
\int_{{\rm domain}\, x}{\rm d}x\,
  p(x-\xi^{\rm ML})\ln p(x-\xi )\, \le\, \int_{{\rm domain}\, x} p(x-\xi^{\rm ML})\ln p(x-\xi^{\rm ML})\, .
    \label{A.7}
\end{equation}

\section{The ML Estimator of the Chi-Squared Model}
\label{B}

The ML estimator of the chi-squared model (\ref{6.11}) is calculated.

Up to an additive constant (independent of $\xi$) the logarithmic likelihood
function is given by
\begin{equation}
  \ln{\cal L}_N(\xi )\, =\, {\rm const}+[x-\xi-\exp (x-\xi )]\, .
  \label{B.1}
\end{equation}
The ML estimator $\xi^{\rm ML}$ solves the ML equation
\begin{eqnarray}
  0\, &=&\, \frac{\rm d}{{\rm d}\xi}\ln{\cal L}_N(\xi )\nonumber\\
  &=&\, -1+\exp (x-\xi )
  \label{B.2}
\end{eqnarray}
The solution is 
\begin{equation}
  \xi^{\rm ML}\, =\, x\, .
  \label{B.4}
\end{equation}

\section{Two Versions of the {F}isher Information}
\label{C}

It shall be shown that the two lines of Eq. (\ref{4.1}) agree with each other.
Let us start from the second line which we write as
\begin{equation}
  \frac{\partial^2}{\partial\xi^2}p(\bold{x}|\xi )\, =\,
  \left.\frac{\partial^2}{\partial\xi^{\prime 2}}
    \int_{-\infty}^{\infty}{\rm d}x\, p(x|\xi )\ln p(x|\xi^{\prime})
    \right|_{\xi^{\prime}=\xi}\, .
    \label{C.1}
\end{equation}
This expression can be rewritten
\begin{eqnarray}
  F(\xi )
  &=&\, -\left.\frac{\partial}{\partial\xi^{\prime}}
  \int{\rm d}x\, p(x|\xi )\, \frac{
    \frac{\partial}{\partial\xi^{\prime}}p(x|\xi^{\prime})}
      {p(x|\xi^{\prime})}\right|_{\xi^{\prime}=\xi}\nonumber\\
      &=&\, -\int{\rm d}x\, p(x|\xi )\left[\frac{
          \frac{\partial^2}{\partial\xi^{\prime 2}}p(x|\xi^{\prime})}
        {p(x|\xi^{\prime})}
        -\left(\frac{
          \frac{\partial}{\partial\xi^{\prime}}p(x|\xi^{\prime})}{p(x|\xi^{\prime})}
        \right)^2\right]_{\xi^{\prime}=\xi}\nonumber\\
        &=&\, -\int{\rm d}x\, \left[\frac{\partial^2}{\partial\xi^2}p(x|\xi )
        -p(x|\xi )\left(\frac{\partial}{\partial\xi}\ln p(x|\xi )\right)^2
        \right]\nonumber\\
      &=&-\frac{\partial^2}{\partial\xi^2}
      \int_{-\infty}^{\infty}{\rm d}x\, p(x|\xi )\, 
      +\, \int_{-\infty}^{\infty}{\rm d}x\,
      p(x|\xi )\left(\frac{\partial}{\partial\xi}
      \ln p(x|\xi )\right)^2\, .
      \label{C.2}
\end{eqnarray}
The first one of the two integrals in the last line vanishes since $p$ is
normalised to unity for every $\xi\, .$ The second integral in the last line corresponds to the
first line of Eq. (\ref{4.1}).

\section{The Chi-Squared Model}
\label{D}

Each of the quantities $x_k\, ,$ where $k=1,2\, ,$ shall have the Gaussian distribution 
\begin{equation}
  w(x_k|\sigma )\, =\, (2\pi\sigma^2)^{-1/2}
  \exp\Big(-\frac{x_k^2}{2\sigma^2}\Big)\, ,\quad\quad -\infty <x_k<\infty\, ,
  \label{D.1}
\end{equation}
with one and the same root mean square value $\sigma\, .$ The chi-squared model $\chi_2^{\rm sq}(T|\sigma )$
with two degrees of freedom is the distribution of the quantity
\begin{equation}
  T\, =\, x_1^2+x_2^2\, .
  \label{D.2}
\end{equation}
It is given by
\begin{equation}
  \chi_2^{\rm sq}(T|\sigma )\, =\, \frac{1}{2\sigma^2\Gamma (1)}
  \exp\Big(-\frac{T}{2\sigma^2}\Big)\, ,\quad\quad 0<T,\tau <\infty\, ,
  \label{D.3}
\end{equation}
see e.g. Eq. (4.34) of Ref. \cite{Harney:16}. This distribution is
normalised to unity. The transformations
\begin{eqnarray}
  z\, &=&\, \ln T\, ,\nonumber\\
  \zeta\, &=&\, \ln (2\sigma^2)
  \label{D.4}
\end{eqnarray}
lead to 
\begin{eqnarray}
  \tilde{\chi}_2^{\rm sq}(z|\zeta )\, &=&\, \frac{{\rm d}T}{{\rm d}z}\, \chi_2^{\rm sq}(T|\sigma )\nonumber\\
  &=&\, T\, \chi_2^{\rm sq}(T|\sigma )\, ,
  \label{D.4a}
\end{eqnarray}
where $T$ and $\sigma$ must be expressed by $z$ and $\zeta\, .$ This gives
\begin{equation}
  \tilde{\chi}_2^{\rm sq}(z|\zeta )\, =\, \, \exp\Big(z-\zeta -e^{z-\zeta}\Big)
  \label{D.5}
\end{equation}
which corresponds to Eq. (\ref{6.11}).

\section{Derivatives of the Functional $H$ for the chi-Squared Distribution}
\label{E}
The integral in the first line of Eq. (\ref{6.12cc}) yields $\Gamma (2)\, .$
It is a special case of the formula
\begin{equation}
  \Gamma (z)\, =\, \int_{-\infty}^{\infty}{\rm d}t\, \exp (zt-e^t)
  \label{E.1}
\end{equation}
given in Sect. 8.312, no. 10 of Ref. \cite{Gradshteyn:15}.
The value of the Gamma function required in Eq. (\ref{6.12cc}) is
\begin{equation}
  \Gamma (2)\, =\, 2\, .
  \label{E.2}
\end{equation}

\section{The Functional $H$ for the Trigonometric Model with Translational Invariance}
\label{F}

The functional $H(\xi^{\rm ML}|\xi )$ for the trigonometric model model (\ref{10.10}) is given by Eq. (\ref{10.12a}) to be
\begin{equation}
  H(\xi^{\rm ML}|\xi )\, =\, {\cal N}\int_{-\pi /2}^{\pi /2}{\rm d}s\, \cos^2(s-\xi^{\rm ML})\ln\frac{\cos^2(s-\xi )}{\cos^2(s-\xi^{\rm ML})}\, .
  \label{F.1}
\end{equation}
The substitution
\begin{equation}
  s^{\prime}\, =\, s-\xi^{\rm ML}
  \label{F.2}
\end{equation}
yields
\begin{eqnarray}
  H(\xi^{\rm ML}|\xi )\, &=&\, {\cal N}\int_{-\pi /2+\xi^{\rm ML}}^{\pi /2+\xi^{\rm ML}}{\rm d}s^{\prime}\, \cos^2s^{\prime}\ln\frac{\cos^2(s^{\prime}-\xi +\xi^{\rm ML})}{\cos^2s^{\prime}}\nonumber\\
    &=&\, {\cal N}\int_{-\pi /2}^{\pi /2}{\rm d}s^{\prime}\, \cos^2s^{\prime}\ln\frac{\cos^2(s^{\prime}-\xi +\xi^{\rm ML})}{\cos^2s^{\prime}}\nonumber\\
    &=&\, H(0|\xi -\xi^{\rm ML})\, .
    \label{F.3}
\end{eqnarray}
Here, the second line is obtained from the first one because the integrand is periodic with a period of $\pi\, ,$
hence, the shift of the limits of integration is immaterial. With the abbreviation
\begin{equation}
  \alpha =\xi^{\rm ML}-\xi
  \label{F.4}
\end{equation}
this reads
\begin{eqnarray}
  H(\xi^{\rm ML}|\xi )\, &=&\, H(0|-\alpha )\nonumber\\
  &=&\, {\cal N}\int_{-\pi /2}^{\pi /2}{\rm d}s\, \cos^2s\ln\frac{\cos^2(s+\alpha)}{\cos^2s}\, .
  \label{F.5}
  \end{eqnarray}
Substituting
\begin{equation}
  s^{\prime}\, =\, -s
  \label{F.6}
\end{equation}
in the integral (\ref{F.5}) one obtains
\begin{eqnarray}
  H(0|-\alpha )\, &=&\, -{\cal N}\int_{\pi /2}^{-\pi /2}{\rm d}s^{\prime}\, \cos^2s^{\prime}\ln\frac{\cos^2(-s^{\prime}+\alpha)}{\cos^2s^{\prime}}\nonumber\\
  &=&\, {\cal N}\int_{-\pi /2}^{\pi /2}{\rm d}s^{\prime}\, \cos^2s^{\prime}\ln\frac{\cos^2(s^{\prime}-\alpha)}{\cos^2s^{\prime}}\nonumber\\
  &=&\, H(0|\alpha )\, .
  \label{F.7}
\end{eqnarray}
Although $\cos(s-\alpha )$ vanishes at a point within the domain of integration, the integral (\ref{F.7}) exists and can be obtained
as if the integrand were simply undefined at this isolated point, see appendix \ref{H}.

Comparing (\ref{F.5}) with (\ref{F.7}) shows that $H(0|-\alpha )$ is a mirror-symmetrical function of the difference $\alpha\, .$
By consequence, all odd derivatives vanish at $\alpha =0\, ,$
\begin{equation}
  \left.\frac{\partial^{2\nu +1}}{\partial\alpha^{2\nu +1}}H(0|-\alpha )\right|_{\alpha =0}\, =\, 0\, ,\quad\quad\quad \nu =0,1,2,\dots\, .
  \label{F.8}
\end{equation}
We calculate the derivatives with respect to $\alpha\, .$ The first derivative is the basis of all higher ones. It must be rewritten
in order to see that all derivatives exist. Starting from Eq. (\ref{F.7}) we find
\begin{eqnarray}
  \frac{\partial}{\partial\alpha}H(0|-\alpha )\, &=&\, {\cal N}\frac{\partial}{\partial\alpha}
  \int_{-\pi /2}^{\pi /2}{\rm d}s\, \left[\cos^2s\ln\cos^2(s-\alpha)\, -\cos^2s\ln\cos^2s\right]\nonumber\\
  &=&\, {\cal N}\frac{\partial}{\partial\alpha}\int_{-\pi /2}^{\pi /2}{\rm d}s\, \cos^2s\ln\cos^2(s-\alpha)\, .
  \label{F.9}
\end{eqnarray}
By use of the substitution
\begin{equation}
  s^{\prime}\, =\, s-\alpha
  \label{F.10}
\end{equation}
one obtains from the integral (\ref{F.9})
\begin{eqnarray}
  \frac{\partial}{\partial\alpha}H(0|-\alpha )\, &=&\, {\cal N}\frac{\partial}{\partial\alpha}
  \int_{-\pi /2}^{\pi /2}{\rm d}s\, \cos^2(s+\alpha )\ln\cos^2s\nonumber\\
    &=&\, 2{\cal N}\frac{\partial}{\partial\alpha}
  \int_{-\pi /2}^{\pi /2}{\rm d}s\, \cos^2(s+\alpha )\ln\cos s\, .
  \label{F.11}
\end{eqnarray}
We express $\cos (s+\alpha )$ by the sum
\begin{equation}
  \cos (s+\alpha )\, =\, \cos s\cos\alpha-\sin s\sin\alpha
  \label{F.12}
\end{equation}
and obtain
\begin{equation}
 \frac{\partial}{\partial\alpha}H(0|-\alpha )\, =\, 2{\cal N}\frac{\partial}{\partial\alpha}
 \int_{-\pi /2}^{\pi /2}{\rm d}s\, \left[\cos s\cos\alpha-\sin s\sin\alpha\right]^2\ln\cos s\, .
 \label{F.13}
\end{equation}
The square $\left[\dots\right]^2$ of a binomial expression displays two squares and a mixed term.
Here, the mixed term, as a function of $s\, ,$ is antisymmetric with respect to $s=0\, .$ Therefore
the integral over the mixed term vanishes and we obtain
\begin{equation}
  \frac{\partial}{\partial\alpha}H(0|-\alpha )\, =\, 2{\cal N}\frac{\partial}{\partial\alpha}
  \int_{-\pi /2}^{\pi /2}{\rm d}s\, \left[\cos^2s\cos^2\alpha+\sin^2s\sin^2\alpha\right]\ln\cos s\, .
 \label{F.14}
\end{equation}
This integrand, as a function of $s\, ,$ is symmetric with respect to $s=0\, .$ Therefore we have
\begin{eqnarray}
\frac{\partial}{\partial\alpha}H(0|-\alpha )\, &=&\, 4{\cal N}\frac{\partial}{\partial\alpha}\cos^2\alpha\int_0^{\pi /2}{\rm d}s\, \cos^2s\ln\cos s
+\, 4{\cal N}\frac{\partial}{\partial\alpha}\sin^2\alpha\int_0^{\pi /2}{\rm d}s\, \sin^2s\ln\cos s\nonumber\\
&=&\, 8{\cal N}\cos\alpha\sin\alpha\left[-\int_0^{\pi /2}{\rm d}s\, \cos^2s\ln\cos s + \int_0^{\pi /2}{\rm d}s\, \sin^2s\ln\cos s\right]\, .
\label{F.15}
\end{eqnarray}

Let us introduce the integrals
\begin{equation}
  I_0\, =\, \int_0^{\pi /2}{\rm d}s\, \ln\cos s
  \label{F.16}
\end{equation}
and
\begin{equation}
  I_2\, =\, \int_0^{\pi /2}{\rm d}s\, \sin^2s\ln\cos s\, .
  \label{F.17}
\end{equation}
They allow to write Eq. (\ref{F.15}) as
\begin{equation}
  \frac{\partial}{\partial\alpha}H(0|-\alpha)\, =\, 8{\cal N}[-I_0+2I_2]\cos\alpha\sin\alpha
  \label{F.18}
\end{equation}
which gives
\begin{equation}
  \frac{\partial}{\partial\alpha}H(0|-\alpha)\, =\, -4{\cal N}[I_0-2I_2]\sin (2\alpha )
  \label{F.19}
\end{equation}
by help of the identity
\begin{equation}
  \sin (2\alpha )\, =\, 2\cos\alpha\sin\alpha\, .
  \label{F.20}
\end{equation}
The derivative (\ref{F.19}) vanishes at $\alpha =0\, ,$ as expected from Eq. (\ref{F.8}).

The values of the integrals (\ref{F.16}) and (\ref{F.17}) can be taken
from the table of integrals \cite{Gradshteyn:15}. According to the entries
4.387 no. 3 and 8.366 no. 1,2 of \cite{Gradshteyn:15}, the integral (\ref{F.16})
has the value
\begin{equation}
  I_0\, =\, -\frac{\pi}{2}\ln 2
  \label{F.21}
\end{equation}
while the entries 4.387 no. 8 and  8.365 no. 1 as well as 8.366 no.1 yield
\begin{equation}
  I_2\, =\, -\frac{\pi}{8}[2\ln 2+1]\, .
  \label{F.22}
\end{equation}
These two values lead to
\begin{equation}
  I_0-2I_2\, =\, \frac{\pi}{4}\, .
  \label{F.23}
\end{equation}

Together with the value of ${\cal N}$ in Eq. (\ref{10.14b}) the derivative (\ref{F.19}) becomes
\begin{equation}
  \frac{\partial}{\partial\alpha}H(0|-\alpha)\, =\, -2\sin (2\alpha )\, .
  \label{F.24}
\end{equation}
Therefore the second derivative of $H(0|-\alpha )$ becomes
\begin{equation}
  \frac{\partial^2}{\partial\alpha^2}H(0|-\alpha)\, =\, -4\cos(2\alpha )
  \label{F.25}
\end{equation}
and the fourth derivative becomes
\begin{equation}
  \frac{\partial^4}{\partial\alpha^4}H(0|-\alpha)\, =\, 16\cos (2\alpha )\, .
  \label{F.26}
\end{equation}
The maximum of the absolute value of the fourth derivative is found at $\alpha =0$ which means
\begin{equation}
  |H^{(4)}|_{\rm max}\, =\, 16\, .
  \label{F.27}
\end{equation}

Note that the second derivative in Eq. (\ref{F.25}) for $\alpha =0$ gives the value
\begin{equation}
  \left.\frac{\partial^2}{\partial\alpha^2}H(0|\alpha)\right|_{\alpha =0}\, =\, -4
  \label{F.28}
\end{equation}
which is, up to its sign, the Fisher information (\ref{10.26}) of the trigonometric model with translational invariance.

\section{The Normalisation of the Trigonometric Model with Translational Invariance}
\label{G}
By partial integration of the $\cos^2$-function one finds
\begin{equation}
  \int_{-\pi  /2}^{\pi /2}{\rm d}s\, \cos^2s\, =\, \big[\sin s\cos s\big]_{-\pi /2}^{\pi /2}\,  +\,
  \int _{-\pi  /2}^{\pi /2}{\rm d}s\, \sin^2s\, .
\label{G.1}
\end{equation}
This leads to
\begin{equation}
  \int_{-\pi  /2}^{\pi /2}{\rm d}s\, \cos^2s\, =\, \int_{-\pi  /2}^{\pi /2}{\rm d}s\, (1-\cos^2s)
  \label{G.2}
\end{equation}
or
\begin{equation}
  2\int_{-\pi  /2}^{\pi /2}{\rm d}s\, \cos^2s\, =\, \pi
  \label{G.3}
\end{equation}
which proves Eq. (\ref{10.14b}).

\section{Integrating over a Logarithmic Divergence}
\label{H}
In the $\epsilon$-interval
\begin{equation}
  -\epsilon\, <\, s^{\prime}-s_0^{\prime}\, <\, \epsilon\, ,\quad\quad \epsilon >0\, ,
  \label{H.1}
\end{equation}
around the point $s_0^{\prime}\, ,$ where $\cos^2$ vanishes, the $\cos^2$-function
behaves as
\begin{equation}
  \cos^2s^{\prime}\, =\, c(s^{\prime}-s_0^{\prime})^2
  \label{H.2}
\end{equation}
since the $\cos^2$-function does not become negative and $\epsilon$ can be chosen
arbitrarily small. Here, $c$ is a positive number. We show that the integral over the
$\epsilon$-interval
\begin{equation}
  \int_{s_0^{\prime}-\epsilon}^{s_0^{\prime}+\epsilon}{\rm d}s^{\prime}\, \ln\cos^2s^{\prime}\, 
  =\,  2\epsilon\ln c\, +\, \int_{s_0^{\prime}-\epsilon}^{s_0^{\prime}+\epsilon}{\rm d}s^{\prime}\, \ln (s^{\prime}-s_0^{\prime})^2
  \label{H.3}
\end{equation}
exists.

For this we rewrite
\begin{eqnarray}
  \int_{s_0^{\prime}-\epsilon}^{s_0^{\prime}+\epsilon}{\rm d}s^{\prime}\, \ln (s^{\prime}-s_0^{\prime})^2
  &=&\, \int_{s_0^{\prime}-\epsilon}^{s_0^{\prime}}{\rm d}s^{\prime}\, \ln (s^{\prime}-s_0^{\prime})^2\, +\,
  \int_{s_0^{\prime}}^{s_0^{\prime}+\epsilon}{\rm d}s^{\prime}\, \ln (s^{\prime}-s_0^{\prime})^2\nonumber\\
  &=&\, 2\int_{s_0^{\prime}-\epsilon}^{s_0^{\prime}}{\rm d}s^{\prime}\, \ln (s_0^{\prime}-s^{\prime})\, +\, 
  2\int_{s_0^{\prime}}^{s_0^{\prime}+\epsilon}{\rm d}s^{\prime}\, \ln (s^{\prime}-s_0^{\prime})\nonumber\\
  \label{H.4}
\end{eqnarray}
so that the arguments of the logarithms are non-negative. In the first integral of the second line
we substitute
\begin{equation}
  s_0^{\prime}-s^{\prime}\, =\, x
  \label{H.5}
\end{equation}
and obtain
\begin{eqnarray}
  2\int_{s_0^{\prime}-\epsilon}^{s_0^{\prime}}{\rm d}s^{\prime}\, \ln (s_0^{\prime}-s^{\prime})\,
  &=&\, -2\int_{\epsilon}^0{\rm d}x\, \ln x\nonumber\\
  &=&\, 2\int_0^{\epsilon}{\rm d}x\, \ln x\nonumber\\
  &=&\, 2\big[x\ln x - 1\big]_{x=0}^{x=\epsilon}\nonumber\\
  &=&\, 2\epsilon\ln\epsilon\, .
  \label{H.6}
\end{eqnarray}
In a similar way one finds the same result for the second integral on the second line of
(\ref{H.4}). Thus Eq. (\ref{H.3}) yields
\begin{equation}
  \int_{s_0^{\prime}-\epsilon}^{s_0^{\prime}+\epsilon}{\rm d}s^{\prime}\, \ln\cos^2s_0^{\prime}\,
  =\, 2\epsilon\ln c\, +\, 4\epsilon\ln\epsilon\, ,
  \label{H.7}
\end{equation}
and this contributes a negligible amount to the expression (\ref{10.11}) when $\epsilon$
is small. Thus $\ln{\cal L}_N$ and the functional $H(\xi^{\rm ML}|\xi )$ exist.

\section{The Likelihood Function of a Gaussian Model}
\label{I}

The $N$-fold Gaussian model
\begin{equation}
  G_N(\bold{x}|\xi )\, =\, (2\pi\sigma^2)^{-N/2}\prod_{k=1}^N\, \exp\left(-\frac{(x_k-\xi )^2}{2\sigma^2}\right)
  \label{I.1}
\end{equation}
is given in Eq. (\ref{1.6b}). We write it as
\begin{eqnarray}
  G_N(\bold{x}|\xi )\, &=&\, (2\pi\sigma^2)^{-N/2}\exp\left(-\frac{1}{2\sigma^2}\sum_{k=1}^N\, (x_k^2-2x_k\xi +\xi^2)\right)\nonumber\\
  &=&\, (2\pi\sigma^2)^{-N/2}\exp\left(-\frac{1}{2\sigma^2}\sum_{k=1}^N\, (<x^2>-2<x>\xi +\xi^2)\right)
  \label{I.2}
\end{eqnarray}
by introducing the averages
\begin{eqnarray}
  <x^2>\, &=&\, \frac{1}{N}\sum_{k=1}^N\, x_k^2\, ,\nonumber\\
  <x>\, &=&\, \frac{1}{N}\sum_{k=1}^N\, x_k\, .
  \label{I.3}
\end{eqnarray}
It is not difficult to factorise this according to
\begin{equation}
  G_N(\bold{x}|\xi )\, =\, (2\pi\sigma^2)^{-N/2}\left(\exp (-\frac{N}{2\sigma^2}(<x^2>-<x>^2)\right)\exp\left(-\frac{N}{2\sigma^2}(<x>-\xi )^2\right)\, .
  \label{I.4}
\end{equation}

The posterior distribution $G_N(\xi |\bold{x})$ is given by the factor that depends on $\xi\, ,$ i.e.
\begin{equation}
  G_N(\xi |\bold{x})\, \propto\, \exp\left(-\frac{N}{2\sigma^2}(<x>-\xi )^2\right)
  \label{I.5}
\end{equation}
in agreement with Eq. (\ref{1.6}). The maximum of the likelihood function occurs at
\begin{equation}
\xi^{\rm ML}\, =\, <x>\, .
\label{I.6}
\end{equation}

\section{The Fisher Information of the Binomial Model}
\label{J}

The Fisher information of the model binomial model (\ref{7.1}) shall be independent of $\xi\, .$
This means that the expression
\begin{eqnarray}
  F(\xi )\, &=&\, \sum_{x=0}^1\, q(x|\xi )\left[\frac{\partial}{\partial\xi}\ln q(x|\xi )\right]^2\nonumber\\
  &=&\, [1-R]\left[\frac{R^{\prime}}{1-R}\right]^2\, +\, R\left[\frac{R^{\prime}}{R}\right]^2 \nonumber\\
  &=&\, \frac{[R^{\prime}]^2}{1-R}\, +\, \frac{[R^{\prime}]^2}{R}\nonumber\\
  &=&\, \frac{[R^{\prime}]^2}{R[1-R]}
    \label{J.1}
\end{eqnarray}
be independent of $\xi\, .$ Here, $R^{\prime}$ is the derivative of $R(\xi )\, .$ 
Thus the numerator of (\ref{J.1}) should be proportional to the denominator. This is reached
when we set
\begin{equation}
  R(\xi )\, =\, \cos^2\xi\, ,\quad\quad\quad  -\pi /2\le \xi\le \pi /2
  \label{J.2}
\end{equation}
and obtain
\begin{equation}
  F(\xi )\, \equiv\, 4\, .
  \label{J.3}
\end{equation}

\end{appendix}
%\bibliography{Bayes,Statistik,Sonstiges,Gruppentheorie}
%\bibliographystyle{plain}

\end{document}